\documentclass[pdflatex,sn-mathphys-num]{sn-jnl}

\usepackage{graphicx}\usepackage{multirow}\usepackage{amsmath,amssymb,amsfonts}\usepackage{mathtools}\usepackage{enumitem}\usepackage{amsthm}\usepackage{mathrsfs}\usepackage[title]{appendix}\usepackage{xcolor}\usepackage{textcomp}\usepackage{manyfoot}\usepackage{booktabs}\usepackage{algorithm}\usepackage{algorithmicx}\usepackage{algpseudocode}\usepackage{listings}\usepackage{comment}
\usepackage{tcolorbox}
\usepackage[ draft,
    commandnameprefix=ifneeded, 
    xcolor=dvipdf,
    todonotes={colorinlistoftodos,
                prependcaption,
                textsize=footnotesize,
                backgroundcolor=orange!10,
                textcolor=black,
                linecolor=orange,
                bordercolor=orange}
]{changes}

\theoremstyle{thmstyleone}\newtheorem{theorem}{Theorem}\newtheorem{proposition}[theorem]{Proposition}\newtheorem{assumption}{Assumption}
\newtheorem{lemma}{Lemma}

\theoremstyle{thmstyletwo}\newtheorem{example}{Example}\newtheorem{remark}{Remark}
\theoremstyle{thmstylethree}\newtheorem{definition}{Definition}
\newif\ifinappendix
\inappendixfalse

\newtheorem{corollary}{Corollary}
\newtheorem{standingAssumption}{Standing Assumption}

\definecolor{mygreen}{RGB}{0,128,0}

\AtEndEnvironment{remark}{\unskip\nobreak\hfill$\blacklozenge$}
\AtEndEnvironment{example}{\unskip\nobreak\hfill$\blacktriangle$}
\definecolor{mymagenta}{RGB}{255,0,255}

\usepackage{import}

\newcommand{\diffFunc}[1]{\mathrm{d}{#1}}

\newcommand{\dom}{\mathrm{dom \:}}

\newcommand{\SOthree}{{\operatorname{SO}(3)}}

\newcommand{\R}[1]{\mathbb{R}^{#1}}

\newcommand{\T}[2]{\mathrm{T}_{#1}{#2}}

\newcommand{\M}{{\mathcal{M}}}
\newcommand{\N}{{\mathcal{N}}}

\newcommand{\W}{\mathcal{W}}

\newcommand{\A}{\mathcal{A}}
\newcommand{\U}{\mathcal{U}}
\newcommand{\V}{\mathcal{V}}
\newcommand{\hybrid}{\mathcal{H}}

\renewcommand{\T}[2]{\mathrm{T}_{#1}{#2}}
\newcommand{\Tcone}[2]{\mathrm{T}^{\M}_{#1}({#2})}
\newcommand{\TconeEuclidean}[2]{\mathrm{T}_{#1}({#2})}

\newcommand{\setto}{\rightrightarrows}
\renewcommand{\dom}[1]{\operatorname{dom}{#1}}
\newcommand{\interior}[1]{\operatorname{int}#1}

\newcommand{\vlift}[2]{\operatorname{vlft}_{#1}(#2)}
\newcommand{\grad}[2]{\operatorname{grad} {#1}({#2})}

\newcommand{\rge}{\operatorname{rge}}
\newcommand{\graph}[1]{\operatorname{gph} {#1}}

\renewcommand{\cal}[1]{\mathcal{#1}}

\makeatletter
\newcommand{\myitem}[1]{%
\item[#1]\protected@edef\@currentlabel{#1}%
}
\makeatother

\begin{document}
 
\title[Article Title]{Multiple Control Barrier Functions for Geometric Hybrid Systems: Modeling and Safe Control Synthesis}

\author*[1]{\fnm{Piyush P.} \sur{Jirwankar}}\email{pjirwank@ucsc.edu}

\author[1]{\fnm{Daniel E.} \sur{Ochoa}}\email{dochoatamayo@ucsc.edu}

\author[1]{\fnm{Carlos A.} \sur{Montenegro G.}}\email{camonten@ucsc.edu}

\author[1]{\fnm{Ricardo} \sur{Sanfelice}}\email{ricardo@ucsc.edu}

\affil*[1]{\orgdiv{Department of Electrical and Computer Engineering}, \orgname{University of California Santa Cruz}, \orgaddress{\city{Santa Cruz}, \postcode{95064}, \state{CA}, \country{USA}}}

\abstract{
    We formulate geometric hybrid dynamical systems with states on intrinsic manifolds and inputs in topological spaces using \emph{relaxed $\text{C}^{\text{\textup{1}}}$-control bundles}. Within this framework, we establish a barrier certificate that guarantees forward invariance of a closed set for closed-loop systems. We introduce multiple control barrier functions and derive conditions that render their sublevel sets forward invariant. We then construct pointwise optimal safe control laws, establish regularity conditions for measurable selections, and use these laws to design safety filters for geometric hybrid systems on manifolds. A second-order bouncing-pendulum example illustrates the framework and results.
}


\maketitle

\section{Introduction}
\label{sec:introduction}

Hybrid dynamical systems provide a natural framework for dynamical systems that combine continuous evolution with discrete events. They arise in networked control~\cite{hespanha2007networked}, impulsive systems~\cite{brogliato2016nonsmooth}, and mechanical and aerospace applications, where impacts, switching logic, resets, or communication events interrupt otherwise continuous evolution of states. Among the numerous modeling approaches for closed-loop hybrid systems~\cite{goebel_hybrid_2012, lygeros2003dynamical, haddad2014impulsive}, the hybrid inclusion framework~\cite{goebel_hybrid_2012} offers a systematic foundation for analysis in Euclidean spaces. Its use of set-valued flow and jump dynamics accommodates nonsmooth behavior as well as nonuniqueness caused by uncertainty or external perturbations.

Many mechanical and aerospace systems, however, have states that evolve naturally on curved spaces, or \emph{manifolds}. Applying Euclidean models to these systems often requires coordinates or an embedding into a larger ambient space, which can obscure their intrinsic geometry and introduce artificial constraints. Recent work~\cite{JirwankarSanfelice_CDC2025,JirwankarOchoaSanfelice_CDC2025,jirwankar2026barrier} extends the hybrid inclusion framework to geometric hybrid systems whose state spaces are intrinsic manifolds. The term ``intrinsic'' highlights that the manifold is defined abstractly and that its modeling and analysis do not depend on a Riemannian metric or on an embedding into a Euclidean space. This generality also admits disconnected manifolds, which often arise when modeling hybrid automata with finitely many logic modes~\cite{JirwankarSanfelice_CDC2025}.

Existing results for geometric hybrid systems include Lyapunov conditions for uniform global asymptotic stability, a hybrid invariance principle on manifolds~\cite{JirwankarSanfelice_CDC2025, JirwankarOchoaSanfelice_CDC2025}, and barrier certificates for forward invariance of closed sets~\cite{jirwankar2026barrier}. These results concern closed-loop systems and, therefore, do not provide a framework for designing controllers for geometric hybrid systems with inputs. Such a framework must account not only for the geometry of the state space, but also for the geometry and state dependence of the admissible inputs. In particular, the set of admissible inputs may vary with the state; see~\cite{nijmeijer1982ControlledInvariance}, \cite[Chapter 13.5]{NijmeijerSchaft1990} for a similar discussion for continuous-time systems on manifolds. This feature is fundamental in geometric control because an input applied at one state may belong to a different space from an input applied at another state. During flows, for example, forces or velocities may lie in state-dependent tangent spaces. 

This modeling challenge is closely connected to safety. Safety specifications in cyber-physical systems are commonly expressed by requiring a designated set of states to remain forward invariant~\cite{Ames2017Barrier}. Control barrier functions provide pointwise conditions for enforcing this requirement and have enabled the development of safety filters~\cite{hsu2023safety,choi2020reinforcement}. A safety filter receives an input from a nominal controller and modifies it only when needed to satisfy the barrier constraints. The filtered input is commonly obtained from a pointwise optimization problem, which preserves the nominal behavior as closely as possible while enforcing safety. The theory of pointwise optimal controllers for hybrid systems in Euclidean spaces is developed in~\cite{montenegro2025pointwise}. Extending this idea to geometric hybrid systems requires optimization over state-dependent admissible input sets during both flows and jumps, together with conditions ensuring that the resulting optimal feedback laws possess the regularity needed for well-posed closed-loop dynamics.

To this end, this paper develops a unified framework for modeling and safely controlling geometric hybrid systems with inputs. The main contributions are as follows.

\begin{enumerate}
    \item We introduce a modeling framework for geometric hybrid systems on manifolds with inputs (Section~\ref{sec:GHDS-inputs}). The model is built on \emph{relaxed $C^1$-control bundles}, which allow the admissible input fiber to vary with the state and relax the standard control-bundle formulation from~\cite{PappasSastry2000,Tabuada2005Quotients}.

    \item We establish barrier-function conditions for forward pre-invariance and forward invariance of closed sets for {geometric hybrid closed-loop systems} (Theorem~\ref{thm:BimpliesForward} and Corollary~\ref{cor:forwardInvariance}). These conditions accommodate multiple, possibly nonsmooth barrier functions and treat both continuous evolution and discrete jumps.

    \item We introduce multiple control barrier functions for geometric hybrid plants (Definition~\ref{def:mCBF}) and construct set-valued flow and jump feedback laws from their pointwise constraints. We show that these feedback laws render the prescribed set forward pre-invariant (Theorem~\ref{thm:mCBFForwardPreInvariant}), and we give additional conditions that guarantee the existence and completeness of solutions (Corollary~\ref{cor:mCBFExistenceCompleteness}).

    \item Given a cost function, we provide conditions under which the admissible feedback maps possess measurable pointwise optimal selections (Assumption~\ref{ass:measurableOptimalSelections} and Theorem~\ref{thm:measurableOptimalSelections}). We also establish conditions for outer semicontinuity and local precompactness of the corresponding set-valued optimal feedback laws (Proposition~\ref{prop:osc}).

    \item We apply these results to geometric hybrid safety filters (Section~\ref{sec:safetyFilter}) and illustrate them through a running example involving a second-order pendulum with impacts at a wall (Examples~\ref{ex:pendulum-plant-dynamics}, \ref{ex:pendulum-mCBF-candidate}, \ref{ex:pendulum-mCBF} and~\ref{ex:pendulum-safety-filter}). This geometric hybrid system is controlled by a force during continuous evolution and by the coefficient of restitution at impact.
\end{enumerate}

\emph{Relation with previous work.} Building on our preliminary conference papers~\cite{jirwankar2026barrier,montenegro2025pointwise}, this article presents the complete theoretical development, including the intermediate results and detailed proofs omitted from those earlier versions. The present paper also improves over~\cite{jirwankar2026barrier, montenegro2025pointwise} in the following ways: (i) a modeling approach for geometric hybrid systems with inputs is presented, and a necessary condition for existence of nontrivial solutions is proposed (Lemma~\ref{lem:necessary}); (ii) a variation of Filippov's Lemma for geometric hybrid systems is presented (Lemma~\ref{lem:hybridInput}); (iii) a nonsmooth multiple control barrier function candidate is defined (Definition~\ref{def:mCBF_candidate}) as opposed to the continuously differentiable candidate in~\cite{montenegro2025pointwise}; (iv) conditions for existence of nontrivial solutions are provided (Corollary~\ref{cor:mCBFExistenceCompleteness}) so that the forward pre-invariance certificate does not hold vacuously; (iii) topological conditions for the existence of a Borel-measurable feedback law are introduced (Assumption~\ref{ass:measurableOptimalSelections} and Theorem~\ref{thm:measurableOptimalSelections}); and (iv) the results are applied in the context of geometric safety filters for hybrid systems (Section~\ref{sec:safetyFilter}). To the best of the authors' knowledge, a similar construction of safety filters for hybrid dynamical systems is absent from the literature.

\section{Notation and Preliminaries}
\label{sec:preliminaries}

The set of real, positive real, and nonnegative real numbers is denoted by $\R{}$, $\R{}_{>0}$, and $\R{}_{\geq 0}$, respectively. The set of natural numbers, including $0$, is denoted by $\mathbb{N}$. Given sets $\A$ and $X$ with $\A\subset X$, a function $f: X \to \R{}_{\geq 0}$ is positive definite with respect to $\A$ if $f(x) = 0$ for each $x\in \A$ and $f(x) > 0$ otherwise. 

\subsection{Topology, Differential Geometry, and Measurable Spaces}

\medskip \noindent
\emph{Topology.} Let $(X, \tau)$ be a topological space. Given a set $S\subset X$, its interior, boundary, and closure are denoted $\interior{S}$, $\partial S$, and $\overline{S}$, respectively. Given a set $\A\subset X$, a set $U$ is a \emph{neighborhood of $\A$} if there exists $U'\in \tau$ such that $\A\subset U' \subset U$. If, additionally, the set $U$ is open (resp., compact), then $U$ is an open (resp., compact) neighborhood of $\A$. A set $S \subset X$ is compact if each open cover of $S$ has a finite subcover, and it is precompact if $\overline{S}$ is compact. 

\medskip\noindent
\emph{Differential geometry.} An $n$-dimensional topological manifold $\M$ is a second-countable Hausdorff space that is locally Euclidean of dimension $n$. A coordinate chart of $\M$ at $x$ is a pair $(U, \varphi)$ where $x\in U\subset \M$, with $U$ an open set and $\varphi:U\to \varphi(U)\subset\R{n}$ a homeomorphism. A collection of coordinate charts $\{(U_{\alpha}, \varphi_{\alpha}) \}_{\alpha\in \mathcal{I}}$, where $\mathcal{I}$ denotes some index set, such that $\bigcup_{\alpha}U_\alpha = \M$ is called an \emph{atlas}.
Two charts $(U, \varphi)$ and $(W, \psi)$ of $\M$ are $C^k$-compatible if $U\cap W = \varnothing$ or $\psi\circ\varphi^{-1}: \varphi(U\cap W)\to \psi(U\cap W)$ is a $C^k$-diffeomorphism. A $C^k$-atlas  of $\M$ is one whose charts are $C^k$-compatible. A $C^k$-manifold  $\M$ is a topological manifold endowed with a maximal $C^k$-atlas. The tangent space to~$\M$ at $x\in\M$ is denoted by $\T{x}{\M}$ and the tangent bundle of $\M$ is denoted by $\T{}{\M}$.  For more details on differential manifolds, we refer the reader to \cite{Lee}. \smallbreak 

A map $f:\M\to\N$ between two $C^k$-manifolds is $C^k$ at~$x\in\M$ if there exists a coordinate chart $(U,\varphi)$ at $x\in \M$ and a chart $(W, \psi)$ at $f(x)\in \N$ such that the coordinate representation $\psi\circ f \circ \varphi^{-1}$ is $C^k$ in the usual Euclidean sense. The map $f$ is said to be $C^k$ if it is $C^k$ at each $x\in\M$. The tangent map of $f$ at $x\in\M$, denoted by~$\diffFunc{f_x}:\T{x}{\M}\to\T{f(x)}{\N}$, is defined by $\diffFunc{f_x}(v)\coloneqq  \left.\frac{df}{dt}(\gamma(t)) \right\vert_{t=0}$, where $\mathcal{I}\ni t\mapsto \gamma(t)$, with $\mathcal{I}\subset \R{}$ such that $0\in\mathcal{I}$, is a smooth curve on $\M$ satisfying $\gamma(0)=x$ and $\gamma'(0)=v\in\T{x}{\M}$. 

\medskip\noindent
\emph{Measurable spaces.}
A \emph{$\sigma$-algebra} on a set $X$ is a collection $\mathfrak{A}$ of subsets of $X$ that contains $\varnothing$, is closed under complements, and is closed under countable unions. The pair $(X, \mathfrak{A})$ is called a \emph{measurable} space. Given measurable spaces $(X, \mathfrak{A})$ and $(Y, \mathfrak{C})$, a function $f : X \to Y$ is \emph{$(\mathfrak{A},\mathfrak{C})$-measurable} if, for each $A \in \mathfrak{C}$, its preimage under $f$ satisfies $f^{-1}(A) \in \mathfrak{A}$. 
Given a topological space $(X, \tau)$, the \emph{Borel $\sigma$-algebra} $\mathfrak{B}(X)$ on $X$ is the $\sigma$-algebra generated by the sets in $\tau$. The Lebesgue $\sigma$-algebra on $\R{n}$, denoted by $\mathfrak{L}(\R{n})$, is the $\sigma$-algebra generated by the collection of sets $\tau_{\R{n}} \cup N$, where $\tau_{\R{n}}$ is the usual metric topology on $\R{n}$, namely, the topology induced by the Euclidean distance, and $N$ is the collection of all subsets of $\R{n}$ having Lebesgue measure zero.

A function $f : X \to Y$ is \emph{Borel measurable} if it is $(\mathfrak{B}(X), \mathfrak{B}(Y))$-measurable. A function $f : \R{n}\to Y$ is \emph{Lebesgue measurable} if it is $(\mathfrak{L}(\R{n}), \mathfrak{B}(Y))$-measurable. 

\subsection{Nonsmooth Analysis on Manifolds}

We recall the following preliminaries on geometric nonsmooth analysis from~\cite{JirwankarSanfelice_CDC2025}. First, given a $C^1$-manifold $\M$ and a set $S\subset \M$, we define the contingent cone to $S$ at $x\in S$. 

\begin{definition}[{Tangent cone to subsets of manifolds~\cite[Definition~4]{JirwankarSanfelice_CDC2025}}]
    \label{def:tangentConeManifolds}
    Let $\M$ be an $n$-dimensional $C^1$-manifold and let $S \subset \M$ be nonempty. The tangent cone $\Tcone{S}{x}$ at $x\in S$ is defined as
    \begin{align*}
        \Tcone{S}{x} \coloneqq \left(\diffFunc{\varphi_x}\right)^{-1}(\TconeEuclidean{\varphi(S\cap U)}{\varphi(x)}),
    \end{align*}
    where $(U, \varphi)$ is any coordinate chart at $x$, and $\TconeEuclidean{\varphi(S\cap U)}{\varphi(x)}$ denotes the Euclidean tangent cone to the set $\varphi(S\cap U)\subset \mathbb{R}^n$ at the point $\varphi(x)$.
\end{definition}

Next, we define the notion of a locally Lipschitz function between $C^1$-manifolds. 

\begin{definition}[{Locally Lipschitz function~\cite[Definition~4.1]{jirwankar2026barrier}}] 
    \label{def:Lipschitz}
    Given two $C^1$-manifolds $\cal M$ and $\cal N$, a function $f: \cal M\to \cal N$ is locally Lipschitz if, for each $x \in \cal M$, there exist a coordinate chart $(U,\varphi)$ on $\cal M$ at $x$ and a coordinate chart $(W,\psi)$ on $\cal N$ at $f(x)$ such that the function $\psi\circ f \circ \varphi^{-1} : \varphi(U\cap f^{-1}(W))  \to \psi(W)$ is locally Lipschitz in the Euclidean sense.
\end{definition}

Given a $C^1$-manifold $\M$, a locally Lipschitz function $f : \M \to \R{}$ is differentiable almost everywhere. At points of nondifferentiability of $f$, the following notion of (Clarke) generalized directional derivative measures its maximum rate of change at $x\in \M$ in the direction $v\in\T{x}{\M}$. 

\begin{definition}[Generalized directional derivative]
    \label{def:generalizedDerivative}
    The generalized directional derivative of a locally Lipschitz function $f: \dom{f}\to\R{}$, $\dom{f}\subset \cal M$, at $x\in\cal M$ in the direction $v\in\T{x}{\cal M}$ is denoted by $f^\circ (x, v)$ and defined as
    \begin{align}\label{eq:defGeneralizedDerivative}
        f^\circ (x, v) \coloneqq \limsup_{\substack{y \to x\\  t \searrow 0}} \frac{\widetilde{f}(\varphi(y) + t\diffFunc{\varphi_x}(v)) - \widetilde{f}(\varphi(y))}{t}
    \end{align}
    where $(U, \varphi)$ is a chart at $x$, $\widetilde{f}\coloneqq f \circ\varphi^{-1}$, and recall that $\diffFunc{\varphi_x}: \T{x}{\cal M}\to \mathbb{R}^{\dim{\cal M}}$ denotes the differential of $\varphi$ at $x$.
\end{definition}

\subsection{Set-Valued Analysis on Manifolds}

The following notions are due to~\cite[Ch.~6]{Beer1993}. A set-valued map $F : X \setto Y$ between topological spaces maps each point $x\in X$ to a set $F(x)\subset Y$. The domain, range, and graph of $F$ are defined as $\dom{F}\coloneqq \{x\in X: F(x)\neq \varnothing\}$, $\rge{F}\coloneqq \{y\in Y : \exists x\in X \text{ s.t. } y\in F(x)\}$, and $\graph{F}\coloneqq \{(x,y)\in X\times Y : y \in F(x)\}$, respectively. A selection from $F$ is a function $f : X\to Y$ that satisfies $\graph{f}\subset \graph{F}$. Given a set $O \subset Y$, the \emph{lower inverse} of $O$ is defined by $F^{-1}(O)\coloneqq \{x\in X : F(x)\cap O \neq \varnothing\}$. 

Now, let $X$ and $Y$ be  Hausdorff spaces. The map $F: X \setto Y$ is \emph{outer semicontinuous} at $x_0\in X$ if $\limsup_{x \to x_0} F(x) \subset F(x_0)$, where the $\limsup$ is computed using nets in~$X$; see~\cite{Beer1993} for more details. 
Given a set $S\subset X$, $F$ is \emph{outer semicontinuous relative to $S$} if the map $\widetilde{F} : X \setto Y$, defined as  $\widetilde{F}(x)\coloneqq F(x)$ for each $x\in S$ and empty otherwise, is outer semicontinuous at each $x_0\in S$.
The map $F : X\setto Y$ is upper semicontinuous at $x_0 \in X$ if, for each open set $V\subset Y$ containing $F(x_0)$, there exists an open neighborhood $U$ of $x_0$ such that $F(U)\subset V$. The map $F$ is \emph{upper semicontinuous relative to $S\subset X$} if $\widetilde{F}$, as defined above, is upper semicontinuous at each $x_0\in S$.

The map $F$ is inner semicontinuous at $x_0\in X$ if, for each open set $V\subset Y$ and each $x_0\in F^{-1}(V)$, there exists an open neighborhood $U$ of $x_0$ such that $U \subset F^{-1}(V)$.
The map $F$ is inner semicontinuous relative to $S \subset X$ if, for each $x_0 \in S$, each $y \in F(x_0)$, and each open neighborhood $V$ of $y$, there exists an open neighborhood $U$ of $x_0$ such that $U \cap S \subset F^{-1}(V)$.

The map $F$ is continuous at $x_0 \in X$ if it is inner semicontinuous and outer semicontinuous at $x_0$. Similarly, $F$ is continuous relative to $S$ if it is inner semicontinuous and outer semicontinuous relative to $S$. Finally, the map $F$ is locally precompact at $x_0\in X$ if there exists a neighborhood $U$ of $x_0$ such that $\overline{F(U)}$ is compact. The map $F$ is locally precompact relative to $S$ if $\widetilde{F}$ is locally precompact at each $x_0\in S$.

\subsection{Control Bundle}
\label{sec:controlBundles}

\begin{figure}[h!]
    \centering
    \includegraphics[width = 0.12\textwidth, trim = 0cm 17.5cm 27cm 0cm]{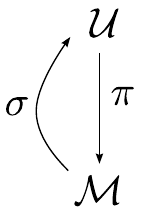}
    \caption{The control bundle $(\M, \U, \pi)$ and a section $\sigma : \M \to \U$. }
    \label{fig:bundle}
\end{figure}

We recall the notion of a \emph{control bundle}, also called a \emph{fibered manifold}, from~\cite[Section~2]{Tabuada2005Quotients} and~\cite{ames2025bundlesbacksteppinggeometriccontrol}. Specifically, a control bundle is a tuple $(\M, \U, \pi)$ of smooth manifolds $\U$ and $\M$, and a surjective submersion\footnote{A function $f : \M \to \N$ between smooth manifolds is a submersion if its differential $\diffFunc{f}_x : \T{x}{\M} \to \T{f(x)}{\N}$ is surjective at each $x\in \M$; see~\cite[Ch.~4]{Lee}. } $\pi : \U \to \M$. For each $x\in \M$, the set $\pi^{-1}(x) \coloneqq \{v \in \U: x = \pi(v)\}$ is called the \emph{fiber over $x$}. A section of the bundle $(\M, \U, \pi)$ is a function $\sigma : \M \to \U$ for which $\pi \circ \sigma = \operatorname{id}_{\M}$, the identity map on $\M$. This construction is depicted in Figure~\ref{fig:bundle}.

\section{Modeling Hybrid Dynamical Systems on Manifolds}

In this section, we present the modeling framework for hybrid dynamical systems on $C^1$-manifolds. In Section~\ref{sec:GHDS}, we recall the framework of \emph{geometric hybrid dynamical systems} from~\cite{JirwankarSanfelice_CDC2025, JirwankarOchoaSanfelice_CDC2025} that models hybrid closed-loop systems on manifolds. We refer to these systems as ``closed geometric hybrid systems''. In Section~\ref{sec:GHDS-inputs}, we extend this framework to model geometric hybrid systems with inputs, which we refer to as ``open geometric hybrid systems.''

\subsection{Geometric Hybrid Dynamical Systems}
\label{sec:GHDS}

We model geometric hybrid dynamical systems using the following inclusion:

\begin{align}
    \label{eq:HS}
    \hybrid : \left\{
        \begin{array}{lccc}
            \dot{x} & \in & F(x) &  \quad x\in C
            \\
            x^+ & \in & G(x) & \quad x\in D
        \end{array}
    \right.
\end{align}
where $x \in \M$ denotes the state of the system and $\M$ is a $C^1$-manifold, $C\subset \M$ and $D\subset \M$ denote the \emph{flow set} and the \emph{jump set}, respectively. The flow map $F : \M \rightrightarrows \T{}{\M}$ and the jump map $G : \M \rightrightarrows \M$ define the set-valued dynamics, where the former is such that $F(x)\subset \T{x}{\M}$ for each $x\in \M$. We denote the data of $\hybrid$ by the notation $\hybrid = (C, F, D, G, \M)$. 

Solutions to $\hybrid$ are defined on hybrid time domains. This notion is defined verbatim from~\cite{goebel_hybrid_2012}. 

\begin{definition}[Hybrid time domain]
    \label{def:HTD}
    A set $E\subset \R{}_{\ge 0}\times \mathbb{N}$ is a compact hybrid time domain if there exists $J\in \mathbb{N}$ such that
    \begin{align}\label{eq:HTD}
    E = \bigcup_{j=0}^{J} [t_j, t_{j+1}]\times  \{j\},
    \end{align}
    for some finite sequence of times $\{t_j\}_{j=0}^{J+1}$ satisfying $0 = t_0 \leq t_1\leq t_2 \leq ... \leq t_{J} \leq t_{J+1}$. A set $E \subset \mathbb{R}_{\geq 0} \times \mathbb{N}$ is a hybrid time domain if it is the union of compact hybrid time domains $E_j$ such that $E_0 \subset E_1 \subset E_2 \subset \ldots \subset E_j \ldots$.
\end{definition}

Next, to define sufficient regularity of the solution during flows, we recall the following notion of local absolute continuity from~\cite{BulloLewis}.

\begin{definition}[{Local absolute continuity~\cite[Sec.~A.2.1]{BulloLewis}}]
    \label{def:localAbsoluteContinuity}
    Let $\M$ be a $C^1$-manifold. A function $f : \R{} \supset I \to \M$ is locally absolutely continuous if, for each $C^1$-function $\gamma : \M \to \R{} $, the composition $\gamma \circ f : I \to \R{}$ is locally absolutely continuous. 
\end{definition}

For a function $f : \R{} \to \R{n}$, local absolute continuity of $f$ is defined using bounds on the distance between function values on compact intervals in $\R{}$. However, as the notion of distance is absent on a $C^1$-manifold $\M$, Definition~\ref{def:localAbsoluteContinuity} bypasses this challenge by requiring the composition $\gamma \circ f$ to be locally absolutely continuous for each $C^1$ test function $\gamma : \M \to \R{}$. 
\begin{definition}[Hybrid arc]
    \label{def:hybridArc}
    A function $\phi : \dom{\phi} \to \M$ is a hybrid arc if $\dom{\phi}$ is a hybrid time domain and, for each $j\in \mathbb{N}$, the function $t\mapsto \phi(t,j)$ is locally absolutely continuous on $I_{\phi}^j\coloneqq \{t : (t,j)\in\dom{\phi}\}$. 
\end{definition}

\begin{definition}[Solution concept]
    \label{def:solution_closedloop}
    A hybrid arc $\phi$ is a solution to $\hybrid = (C, F, D, G, \M)$ if $\phi(0,0) \in \overline{C} \cup D$ and the following conditions hold:
    \begin{enumerate}[label=(S\arabic*)]
        \item for each $j \in \mathbb{N}$ such that $I^{j}_{\phi}$ has a nonempty interior, 
        \begin{align*}
        \begin{array}{rcl}
            \phi(t, j) \in C\phantom{(\phi(t,j))} & \quad \text{for all} \quad & t\in \interior{ I^j}, \\
            \displaystyle{\frac{d\phi}{dt}}(t, j) \in F(\phi(t, j)) & \quad \text{for almost all} \quad & t\in I^j;
        \end{array}
    \end{align*}
        \item for each $(t,j)\in\dom{\phi}$ such that $(t,j+1)\in\dom{\phi}$, 
        \begin{align*}
            \begin{array}{ccc}
                 \phi(t,j)\in D & \text{and} & \phi(t,j+1) \in G(\phi(t,j)). 
            \end{array}
        \end{align*}
    \end{enumerate}
\end{definition}

A solution $\phi$ to $\hybrid$ is \emph{complete} if $\dom{\phi}$ is unbounded. It is maximal if $\phi$ cannot be extended; i.e., there does not exist a solution $\phi'$ to $\hybrid$ with $\dom{\phi'}$ strictly larger than $\dom{\phi}$ such that $\phi'(t,j)=\phi(t,j)$ for each $(t,j)\in\dom{\phi}$. 

Next, we present the following notion of forward invariance of a set $K \subset \M$ for $\hybrid$ that will be used in the forthcoming sections.  

\begin{definition}[Forward pre-Invariance]
    A set $K \subset \M$ is \emph{forward pre-invariant} for $\hybrid$ if each solution $\phi$ to $\hybrid$ with $\phi(0,0)\in K$ satisfies $\rge{\phi}\subset K$. The set $K$ is \emph{forward invariant} for $\hybrid$ if, additionally, each maximal solution to $\hybrid$ starting from $K$ is also complete. 
\end{definition}

In the above definition, the prefix ``pre-'' enables maximal solutions to not be complete.

\subsection{
Modeling Systems with Inputs on Relaxed $C^1$-Control Bundles
}

Extending ideas from~\cite{Tabuada2005Quotients} where a continuous-time system with inputs is modeled on a control bundle, we present the following notion of a relaxed $C^1$-control bundle which captures continuous-time and discrete-time plant dynamics under relaxed assumptions compared to~\cite{Tabuada2005Quotients}.

\begin{definition}[Relaxed $C^1$-control bundle]
    A tuple $(\M, \U, \pi)$ is a relaxed $C^1$-control bundle if $\M$ is a $C^1$-manifold, $\U$ is a topological space, and $\pi : \U \to \M$ is a surjection. 
\end{definition}

\begin{remark}
    In going from a control bundle (see Section~\ref{sec:controlBundles}) to a relaxed $C^1$-control bundle, we do not ask for $\U$ to be a smooth manifold. Consequently, we also lose the ability to ask for differentiability of the map $\pi$. In some of the forthcoming results in this paper, we impose additional structure on $\U$ as needed, thereby relating properties of the dynamical system to regularity of $\U$.  
\end{remark}

Before the introduction of control bundles, continuous-time plants on manifolds were modeled on fiber bundles~\cite{nijmeijer1982ControlledInvariance, NijmeijerSchaft1990}. The following example shows that fiber bundles are also relaxed $C^1$-control bundles.
\begin{example}[Fiber bundles are control bundles]
    \label{ex:fiber-bundle}
    Following~\cite[Chapter 13.5]{NijmeijerSchaft1990}, a fiber bundle is a collection of smooth manifolds $B$, $M$, $U$, a surjective submersion $\pi: B \to M$, and an open cover $\{O_i\}_{i\in \mathcal{I}}$ of $M$ such that, for each $i\in \mathcal{I}$, there exists a diffeomorphism $\Phi_{i} : \pi^{-1}(O_i) \to O_i \times U$ satisfying $\pi(x) = (\operatorname{pr}_1 \circ \Phi_i)(x)$ for each $x\in \pi^{-1}(O_i)$, where $\operatorname{pr}_1 : M \times B \to M$ denotes the projection onto $M$. A fiber bundle is denoted by $(B, M, \pi, U, \{O_i\}_{i\in \mathcal{I}})$.
    Clearly, a fiber bundle is also a relaxed $C^1$-control bundle $(M, B, \pi)$ since $M$ is a $C^1$-manifold and $\pi$ is a surjection, where, in addition, the control input in each fiber belongs to a set that is diffeomorphic to $U$.
\end{example}

\begin{example}[Vector bundle]
    \label{ex:vectorBundle}
    Following~\cite[p.390]{NijmeijerSchaft1990}, a \emph{vector bundle} is a fiber bundle $(B,M,\pi,U,\{O_i\}_{i\in\mathcal{I}})$ where $U$ is a finite-dimensional vector space, the fiber $\pi^{-1}(x)$ is a vector space for each $x\in M$, and, for each $i\in\mathcal{I}$ and each $x\in O_i$, the restriction $\left.\Phi_i\right|_{\pi^{-1}(x)}:\pi^{-1}(x)\to\{x\}\times U$ is a linear isomorphism. As a vector bundle is also a fiber bundle by definition, it is also a relaxed $C^1$-control bundle $(M, B, \pi)$. 
\end{example}

\subsubsection{State-Input Information as Points in the Bundle}
\label{sec:state-input}

Before presenting the modeling of dynamical systems with inputs, we recall from~\cite[Chapter~13]{NijmeijerSchaft1990} how the state-input information is captured on a relaxed $C^1$-control bundle. For a dynamical system with inputs on the relaxed $C^1$-control bundle $(\M,\U,\pi)$, the state belongs to $\M$, whereas a point in $\U$, called a state-input point, encodes the joint information about the state and its associated input.
More precisely, for each state ${z}\in\M$, the set of state-input points is the fiber
\(
    \pi^{-1}({z})\subset\U.
\)
{As each point $\chi \in \U$ belongs to a fiber of $\pi$, its associated state can be obtained using knowledge of $\chi$ only, via
\begin{align}
    \label{eq:state-using-input}
    z=\pi(\chi).
\end{align}
Therefore, the point ${\chi}\in\U$ captures both the state and the input.

Comparing the above approach to dynamical systems with inputs in the Euclidean space~\cite{Sanfelice_HybridFeedbackControl_2021}, the set of admissible control inputs at each state $z\in \R{n}$ is all of $\R{m}$, making input $u\in \R{m}$ independent of the state $z \in \R{n}$. However, such Euclidean modeling can very well be captured in the above control-bundle-based modeling by letting $\M \coloneqq \R{n}$, $\U \coloneqq \M \times \R{m}$, and $\pi : \R{n}\times \R{m} \to \M$ as the projection onto the first component. Then, for each state $z\in \R{n}$, the fiber $\pi^{-1}(z) = \{z\}\times \R{m}$ ensures that the admissible input space is $\R{m}$, and, with $u\in \R{m}$, a point $\chi = (z, u) \in \pi^{-1}(z)$ in the fiber captures the state-input information. The product structure of $\U = \M \times \R{}$ in the Euclidean setting enables representation of a state-input point $\chi \in \U$ by the tuple $(z,u)$.

\begin{remark}[On nomenclature of state-input points]
    Given a relaxed $C^1$-control bundle, we choose to call en element of $\U$ a ``state-input point.'' Doing so encodes the fact that the admissible inputs to the system are state-dependent. The terminology ``point'' in ``state-input point'' also highlights that any $\chi\in \U$ may not be decomposed into a state component and an input component. In general, such a decomposition is possible on fiber bundles where \emph{local trivialization}~\cite[p.268]{Lee} allows locally representing a neighborhood of each point $\chi$ as a product space, thereby enabling the tuple representation of $\chi$ in this neighborhood. Example~\ref{ex:attitude-kinematics} in Appendix~\ref{app:modeling-CT-DT} illustrates one such example, where the local trivialization is, in fact, global. On arbitrary relaxed $C^1$-control bundles, however, such decomposition may not be possible. Therefore, the combined state-input information must be denoted by a point in $\U$.
\end{remark}

In general, such representation is possible on fiber bundles where \emph{local trivialization}~\cite[p.268]{Lee} allows locally representing a neighborhood of each point $\chi$ as a product space, thereby enabling the tuple representation of $\chi$ in this neighborhood. Example~\ref{ex:attitude-kinematics} in Appendix~\ref{app:modeling-CT-DT} illustrates one such example, where the local trivialization is, in fact, global. On arbitrary relaxed $C^1$-manifolds, however, such decomposition of the state-input point into a state component and an input component may not be possible. Therefore, the combined state-input information must be denoted by a point in $\U$.

}

\subsection{Geometric Hybrid Systems with Inputs}
\label{sec:GHDS-inputs}

We model a geometric hybrid system with inputs on a relaxed $C^1$-control bundle $\mathcal{B} = (\M, \U, \pi)$ as
\begin{align}
    \label{eq:hybridplant}
    \hybrid_P : \left\{
        \begin{array}{ll}
           \dot{{z}}\phantom{^+} \in F_P(\chi) & \quad \chi \in C_P
           \\
            {z}^+  \in G_P(\chi) & \quad  \chi \in D_P
            \\
            y\phantom{^+} = {h_P}({\chi}) & 
        \end{array}
    \right.
\end{align}
where ${z}\in \M$ is the state of the plant and ${\chi}\in \U$ denotes the state-input point. This point captures information of the state $z$ using~\eqref{eq:state-using-input} so that $z = \pi(\chi)$; see also Appendix~\ref{app:modeling-CT-DT} for a detailed discussion on such modeling. The flow set $C_P \subset \U$ and the jump set $D_P \subset \U$ capture constraints on the hybrid plant $\hybrid_P$. The flow map $F_P : \U \setto \T{}{\M}$, with $F_P(\chi)\subset \T{\pi(\chi)}{\M}$ for each $\chi\in \U$, captures the continuous evolution (or flow) of the state under the action of $\chi$. Similarly, the jump map $G_P: \U\setto \M$ captures the discrete behavior (or jumps) under the action of $\chi \in D_P$. The output is denoted by $y$ and captured by the output map ${h_P} : \U \to \cal Y$, where $\cal Y$ denotes the output space. The data of $\hybrid_P$ is denoted by $\hybrid_P = (C_P, F_P, D_P, G_P, {h_P}, \cal B)$. The subscript ``$P$'' highlights that $\hybrid_P$ is a \emph{plant}. Recall that we refer to $\hybrid_P$ as an open geometric hybrid system. 

\begin{remark}
    As $\chi\in\U$ encodes state-input information, the maps $F_P$ and $G_P$ do not depend explicitly on the state $z$; see Appendix~\ref{app:modeling-CT-DT} for a detailed discussion.  One could, however, artificially introduce explicit dependence by modeling the flow and jump dynamics differently. In particular, consider the map $\Phi:C_P \to \M\times \U$, defined as $\Phi(\chi) = (\pi(\chi), \chi)$ for each $\chi\in C_P$. Define $\widetilde{C}_P \coloneqq \Phi(C_P)$, and note that $\Phi$ is a bijection from $C_P$ to $\widetilde{C}_P$. Then, with $\widetilde{F}_P : \M \times \U \setto \T{}{\M}$ as $\widetilde{F}_P(z,\chi) \coloneqq F_P(\Phi^{-1}(z,\chi))$ for each $(z,\chi)\in \widetilde{C}_P$ and empty otherwise, the resulting dynamics $\dot{z}\in \widetilde{F}_P(z, \chi)$ $(z, \chi)\in \widetilde{C}_P$ equivalently model the differential inclusion in~\eqref{eq:hybridplant}. That $\Phi$ is bijection enforces both modeling approaches to be equivalent. 
\end{remark}

\subsubsection{Why Relaxed $C^1$-Control Bundles?}

Each ingredient of the relaxed $C^1$-control bundle $\mathcal{B}=(\M,\U,\pi)$ plays a distinct role in modeling $\hybrid_P$. The $C^1$-manifold $\M$ provides the intrinsic state space and its tangent spaces, so that the flow condition $F_P(\chi)\subset\T{\pi(\chi)}{\M}$ and, therefore, the differential inclusion in~\eqref{eq:hybridplant} are well-defined. The topological space $\U$ collects the admissible state-input points, while its topology supplies the notions of continuity of functions, closedness and compactness of sets, and Borel measurability used to impose regularity on the plant data and on feedback laws. The projection $\pi$ associates each ${\chi}\in\U$ with its state ${z}=\pi({\chi})$, and its fiber $\pi^{-1}({z})$ captures the set of state-input points. Finally, surjectivity of $\pi$ guarantees that every state has at least one admissible {state-input point} in $\U$. Thus, the relaxed bundle contains the minimal structure needed to describe the state, its associated state-input space, and the continuous and discrete dynamics of $\hybrid_P$.

The above modeling choice raises a natural question: 
\begin{quote}
    \begin{center}
        \emph{Are there other modeling approaches that are better suited to capture the geometry of $\hybrid_P$?}
    \end{center}
\end{quote}
Existing geometric control literature~\cite{van1982observability,nijmeijer1982ControlledInvariance,Sussmann1976,Brockett2014Survey, NijmeijerSchaft1990} provides two closely related alternatives. Fiber-bundle models use local trivializations to identify the {admissible state-input points} locally with a product space~\cite[Ch.~13]{NijmeijerSchaft1990}, \cite{BARBEROLINAN2020100935}; this is useful when locally uniform fibers are needed for a coordinate-based analysis of the system, but it is stronger than required for the purpose of this paper. Control-bundle models replace local triviality by the assumption that $\U$ is a manifold and $\pi$ is a surjective submersion~\cite{Tabuada2005Quotients}; such smoothness assumptions support differential-geometric comparisons of control systems, whereas the analysis in this paper does not differentiate with respect to the input or the projection. {In both of these modeling approaches, $\U$ is a manifold. However, this assumption can be quite restrictive for hybrid systems, where the inputs during flows and at jumps can have different dimensions and, therefore, may not belong to the same manifold. Allowing $\U$ to be a topological space enables such modeling for hybrid systems.

}

Following the topological modeling approach instead of the control/fiber bundle approach, a third, variational analysis formulation places $({z},u)$ in a product $\M\times\mathcal{V}$, where $\V$ is some topological space, and represents state-dependent input constraints by a set-valued map $U:\M\rightrightarrows\mathcal{V}$~\cite{Sanfelice_HybridFeedbackControl_2021}. In that formulation, the admissible pairs 
form $\graph{U}$ and the admissible inputs at ${z}$ are $U({z})$; equivalently, one may take $\U=\graph{U}$ and $\pi({z},u)={z}$, so that $\pi^{-1}({z})=\{{z}\}\times U({z})$. The relaxed $C^1$-control bundle adopts this intrinsic fiberwise structure directly, without imposing local triviality, a manifold structure on $\U$, or a submersion condition on $\pi$. It therefore captures state-dependent input spaces as fibers over ${z}$ while retaining only the structure needed in the subsequent analysis.

\subsubsection{Solution Concept for $\hybrid_P$}

\begin{figure}[h!]
    \centering
    \includegraphics[width = 0.3\textwidth, trim = 0cm 16cm 22cm 0cm]{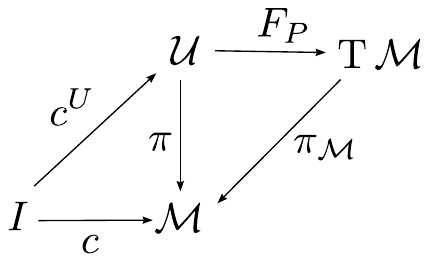}
    \caption{Solution concept from~\cite[Def.~2.3]{Tabuada2005Quotients}. This notion is not well-suited for hybrid systems, where nonsmooth dynamics are prevalent.}
    \label{fig:solutionConcept_Tabuada}
\end{figure}

In this section, we define the notion of a solution to $\hybrid_P$. Note that, when $D_P = \varnothing$ and $F_P$ is single valued, $\hybrid_P$ denotes a continuous-time plant. Such systems on manifolds have been extensively studied in the literature~\cite{Tabuada2005Quotients, ames2025bundlesbacksteppinggeometriccontrol}. Following~\cite[Def.~2.3]{Tabuada2005Quotients}, a solution to such a continuous-time system $\hybrid_P$ is a smooth curve $c : \R{}_{\geq 0} \supset I  \to \M$ for which there exists a (possibly nonsmooth) function $c^U : I \to \U$ such that Figure~\ref{fig:solutionConcept_Tabuada} commutes; namely, $c(t) = (\pi\circ c^{U})(t) = (\pi_{\M}\circ F_P\circ c^{U})(t)$ for each $t\in I$, and $\pi(\chi) = (\pi_{\M}\circ F_{P})(\chi)$ for each $\chi\in \U$. Such solution notion requires smoothness of the state trajectory $c$, and imposes no regularity on the control trajectory $c^U$. This notion of a solution concept is not suitable for hybrid systems applications, where nonsmooth dynamics often force the state trajectory to be nonsmooth. To address this generality in the modeling framework, we present the notion of a solution to~$\hybrid_P$ that extends~\cite[Def.~2.29]{Sanfelice_HybridFeedbackControl_2021} to the geometric setting. 

First, we endow the sets of interest with the respective Borel $\sigma$-algebras.

\begin{standingAssumption}[Borel $\sigma$-algebras on $\M$ and $\U$]
    \label{SA:1}
    Given a relaxed $C^1$-control bundle $(\M, \U, \pi)$, $\M$ and $\U$ are endowed with the respective Borel~$\sigma$-algebras, making $(\M, \mathfrak{B}(\M))$ and $(\U, \mathfrak{B}(\U))$ measurable spaces.
\end{standingAssumption}

The above standing assumption is mild as it enables the definition of a Lebesgue measurable function from a hybrid time domain to $\U$. Such a function, with the following additional properties, seeks to emulate~\cite[Definition~2.27]{Sanfelice_HybridFeedbackControl_2021} to model {hybrid inputs} on $\U$. Recalling from Definition~\ref{def:HTD} and Definition~\ref{def:hybridArc}, respectively, the notions of a hybrid time domain and a hybrid arc, a hybrid input is defined as follows.

\begin{definition}[Hybrid input]
    \label{def:hybridInput}
    A function ${\chi} : \dom{{\chi}} \to \U$ is a hybrid input if $\dom{{\chi}}$ is a hybrid time domain and if, for each $j\in\mathbb{N}$, the function $t \mapsto {\chi}(t,j)$ is Lebesgue measurable and locally essentially precompact\footnote{A function $f: \R{} \to \U$ is locally essentially precompact if, for each $s\in\mathbb{R}$, there exists an open neighborhood $I$ of $s$ and a compact set $K\subset \U$ such that $f(t) \in K$ for almost all $t\in I$. The notions of open neighborhoods and compact sets on $\U$ are well-defined as $\U$ is a topological space by assumption.} on the interval $I^j_{{\chi}} \coloneqq \{t : (t,j)\in\dom{{\chi}}\}$. 
\end{definition}

\begin{definition}[Solution concept]
    \label{def:solution}
Let $\cal B = (\M, \U, \pi)$ be a relaxed $C^1$-control bundle. A pair $(z, \chi)$ of a hybrid arc $z : \dom{z} \to \M$ and a hybrid input ${\chi} : \dom{{\chi}} \to \U$ is a solution to $\hybrid_P = (C_P, F_P, D_P, G_P, {h_P}, \cal B)$ if $\dom{(z,\chi)} \coloneqq \dom{z} = \dom{\chi}$ and the following hold:
    \begin{enumerate}[label =(I\arabic*)]
        \item \label{item:solution-initialCondition} $\chi(0,0)\in \overline{C_P} \cup D_P$;
        \item \label{item:solution-compatibility} $z(t,j) = \pi(\chi(t,j))$ for each $(t,j)\in\dom{(z, \chi)}$;
        \item \label{item:solution-flowCondition} for each $j\in \mathbb{N}$ such that $I^j_{{\chi}}$ has a nonempty interior, ${\chi}$ satisfies
        \begin{align*}
            \begin{array}{ccl}
                \chi(t,j) \in C_P & \text{for all} & t\in \interior{I^j_{{\chi}}},
                \\[7pt]
                \displaystyle{\frac{d z}{dt} (t,j)} \in F_P(\chi(t,j)) & \quad \text{for almost all} \quad  & t\in I^j_{{\chi}};
            \end{array}
        \end{align*}
        \item \label{item:solution-jumpCondition} for each $(t,j)\in \dom{(z, \chi)}$ such that $(t,j+1)\in \dom{(z, \chi)}$, 
        \begin{align*}
            \chi(t,j) \in D_P \qquad \text{and} \qquad z(t,j+1) \in G_P(\chi(t,j)).
        \end{align*}
    \end{enumerate}
\end{definition}

Much like the Euclidean setting~\cite[Lemma~2]{maghenem2022sufficientconditionsforwardinvariance} for hybrid closed-loop systems, each nontrivial solution to the {open geometric hybrid system} $\hybrid_P$ satisfies the following tangent cone condition during flows, {where we recall the geometric notion of tangent cone from Definition~\ref{def:tangentConeManifolds}. }

\begin{lemma}[Necessary condition]
    \label{lem:necessary}
    Let $(z, \chi)$ be a solution to $\hybrid_P = (C_P, F_P, D_P, G_P, {h_P}, \mathcal{B})$, and pick any $j \in\mathbb{N}$ for which $\interior{I^j_{{\chi}}} \neq \varnothing$. Then,
    \begin{align}
        F_P(\chi(t,j)) \cap \Tcone{\pi(C_P)}{z(t,j)} \neq \varnothing \qquad \text{for almost all } t\in I^{j}_{{\chi}}. 
    \end{align}
\end{lemma}

\begin{proof}
    Pick any $j\in\mathbb{N}$ for which $\interior{I^{j}_{\chi}}\neq \varnothing$. Define $y(s)\coloneqq z(s,j)$ for each $s\in I^j_{{\chi}}$. As $z$ is a hybrid arc by Definition~\ref{def:solution}, it is locally absolutely continuous on $I^{j}_{\chi}$. Therefore, $y$ is locally absolutely continuous on $I^{j}_{\chi}$. Consequently, there exists a Lebesgue measure zero set $N\subset I^j_{{\chi}}$ such that, for each $t\in \interior{I^j_{{\chi}}}\setminus N$, $\frac{dy}{dt}(t)$ exists and satisfies
    \begin{align}\label{eq:necessaryConditionVelocity}
        \frac{dy}{dt}(t)\in F_P(\chi(t,j)).
    \end{align}
    Fix any such $t$ and note that $y(t)\in\pi(C_P)$. Pick a coordinate chart $(U,\varphi)$ at $y(t)$. Then, as $y$ is continuous due to being locally absolutely continuous,, there
    exists $\delta>0$ such that $(t-\delta,t+\delta)\subset\interior{I^j_{{\chi}}}$ and $y(s)\in U$ for each $s\in (t-\delta, t+\delta)$.
    Choose any sequence $\{h_k\}_{k\in\mathbb N}\subset(0,\delta)$ such
    that $h_k\to0$. For every $k\in\mathbb N$, define $v_k\coloneqq {\varphi(y(t+h_k))-\varphi(y(t))}/{h_k}$. 
    Then, $\varphi(y(t))+h_kv_k=\varphi(y(t+h_k))$ belongs to
    $\varphi(\pi(C_P)\cap U)$ for each $k\in\mathbb N$. Additionally, as $y$ is differentiable at $t$ and $\varphi$ is a $C^1$-diffeomorphism onto its image,
    \begin{align*}
        \lim_{k\to\infty}v_k
        =\lim_{k\to\infty}\frac{\varphi(y(t+h_k)) - \varphi(y(t))}{h_k}
        =\frac{d(\varphi\circ y)}{dt}(t)
        =\diffFunc{\varphi_{y(t)}}\left(\frac{dy}{dt}(t)\right).
    \end{align*}
    By the definition of the Euclidean tangent cone, we have that
    \begin{align*}
        \diffFunc{\varphi_{y(t)}}\left(\frac{dy}{dt}(t)\right)
        \in
        \TconeEuclidean{\varphi(\pi(C_P)\cap U)}{\varphi(y(t))}.
    \end{align*}
    Consequently, Definition~\ref{def:tangentConeManifolds} yields
    $\frac{dy}{dt}(t)\in\Tcone{\pi(C_P)}{y(t)}$. This inclusion, together
    with~\eqref{eq:necessaryConditionVelocity}, results in
    \begin{align*}
        \frac{dy}{dt}(t)\in
        F_P(\chi(t,j))\cap
        \Tcone{\pi(C_P)}{z(t,j)}.
    \end{align*}
    As the above inclusion holds for each $t\in \interior{I^{j}_{{\chi}}}\setminus N$, and the sets $N$ and $\partial I^{j}_{{\chi}}$ have Lebesgue measure zero, the above inclusion holds for almost all $t\in I^{j}_{{\chi}}$. 
\end{proof}

The following example illustrates the modeling of an open geometric hybrid system on a relaxed $C^1$-control bundle. It models a pendulum rotating on the unit circle $\mathbb{S}^1$. It is an extension of~\cite{jirwankar2026barrier} to open geometric hybrid systems. To this end, following the construction in~\cite{jirwankar2026barrier}, consider the unit circle $\mathbb{S}^1$ as the quotient space $[0,1]/\{0,1\}$ with the equivalence relation $0\sim 1$. We equip $\mathbb{S}^1$ with a Riemannian metric $g_{\mathbb{S}^1}$, causing it to be a smooth Riemannian manifold. Moreover, $\T{}{\mathbb{S}^1}$ is a smooth manifold~\cite[Prop.~3.18]{Lee}. Elements of $\T{}{\mathbb{S}^1}$ are written as $([z_1],[z_2])$, where $[z_1]\in \mathbb{S}^1$, $[z_2]\in \T{[z_1]}{\mathbb{S}^1}$, and $[\cdot]$ denotes the corresponding equivalence class. Let $\pi_{\mathbb{S}^1}  : \T{}{\mathbb{S}^1} \to \mathbb{S}^1$ denote the canonical tangent bundle projection map that, for each $([z_1], [z_2])\in\T{}{\mathbb{S}^1}$, yields $\pi_{\mathbb{S}^1}([z_1], [z_2]) \coloneqq [z_1]$.

\begin{example}[Pendulum with impacts]
    \label{ex:pendulum-plant-dynamics}
    We model a pendulum rotating on the unit circle $\mathbb{S}^1$, that contains a connected wall $\cal W$ of thickness $2\alpha$, $\alpha \in (0, 1/2)$, defined as $\cal W \coloneqq ([0, \alpha] \cup [1-\alpha, 1])/\{0,1\}$. The position $[z_1]\in \mathbb{S}^1$, with $z_1 \in [0,1]$, and the velocity $[z_2]\in \T{[z_1]}{\mathbb{S}^1}$, with $z_2 \in \R{}$, constitute the state $z\coloneqq ([z_1], [z_2])\in \M$ of the pendulum, where $\M \coloneqq \T{}{\mathbb{S}^1}$. This construction of $\mathbb{S}^1$ and $\T{}{\mathbb{S}^1}$ is illustrated in Figure~\ref{fig:pendulum-flowAndJump}. 

    At each position $[z_1]\in \mathbb{S}^1$ when the pendulum is away from the wall $\cal W$, the pendulum evolves continuously under the action of friction force and a control force in $\T{[z_1]}{\mathbb{S}^1}$. At impacts, a control input taking values in $[0,1]$ acts as the coefficient of restitution that scales the velocity after jumps. The constraints are captured in the forthcoming flow and jump sets using the notion of tangent cone in Definition~\ref{def:tangentConeManifolds}.
    The above dynamics are modeled on a relaxed $C^1$-control bundle $(\M, \U, \pi)$. The set $\U$ is defined using the fiber product of $\T{}{\mathbb{S}^1}$ with itself, denoted by the Whitney sum~$\oplus$~\cite[Example~10.7]{Lee}, as
    \begin{align}
        \U &\coloneqq \left(\T{}{\mathbb{S}^1} \oplus \T{}{\mathbb{S}^1}\right)\times[0,1] \nonumber
        \\
        &= \{(z,w,\lambda) \in \T{}{\mathbb{S}^1} \times \T{}{\mathbb{S}^1}\times[0,1] : \pi_{\mathbb{S}^1}(z) = \pi_{\mathbb{S}^1}(w)\}. 
    \end{align}
    As $\U \subset \T{}{\mathbb{S}^1}\times \T{}{\mathbb{S}^1}\times [0,1]$ and leveraging the product structure, for each $\chi \in \U$, there exists a unique $(z,w,\lambda) \in \T{}{\mathbb{S}^1} \times \T{}{\mathbb{S}^1}\times[0,1]$ such that $\chi = (z,w,\lambda)$. For each element $\chi = (z,w, \lambda)\in \U$, its first component $z$ denotes the state of the system, the second component $w$ denotes the control force during flows, and the third component $\lambda$ denotes the control input at jumps. The above construction of $\U$ ensures that, for each $[z_1]\in \mathbb{S}^1$, the velocity and the control force belong to the same tangent space $\T{[z_1]}{\mathbb{S}^1}$. Additionally, $\U$ is a topological space when endowed with the subspace topology of the product topology. Define the projection maps $\operatorname{pr}_1:\U\to\T{}{\mathbb{S}^1}$, $\operatorname{pr}_2:\U\to\T{}{\mathbb{S}^1}$, and $\operatorname{pr}_3:\U\to[0,1]$ as $\operatorname{pr}_1(\chi)=z$, $\operatorname{pr}_2(\chi)=w$, and $\operatorname{pr}_3(\chi)=\lambda$, respectively, for each $\chi=(z,w,\lambda)\in\U$. We define $\pi : \U \to \M$ by $\pi(\chi)=\operatorname{pr}_1(\chi)$. Since this projection is surjective, $(\M, \U, \pi)$ is a relaxed $C^1$-control bundle.

    \begin{figure}
        \centering
        {\includegraphics[width=0.8\linewidth, trim = 0cm 5cm 1cm 0cm, clip]{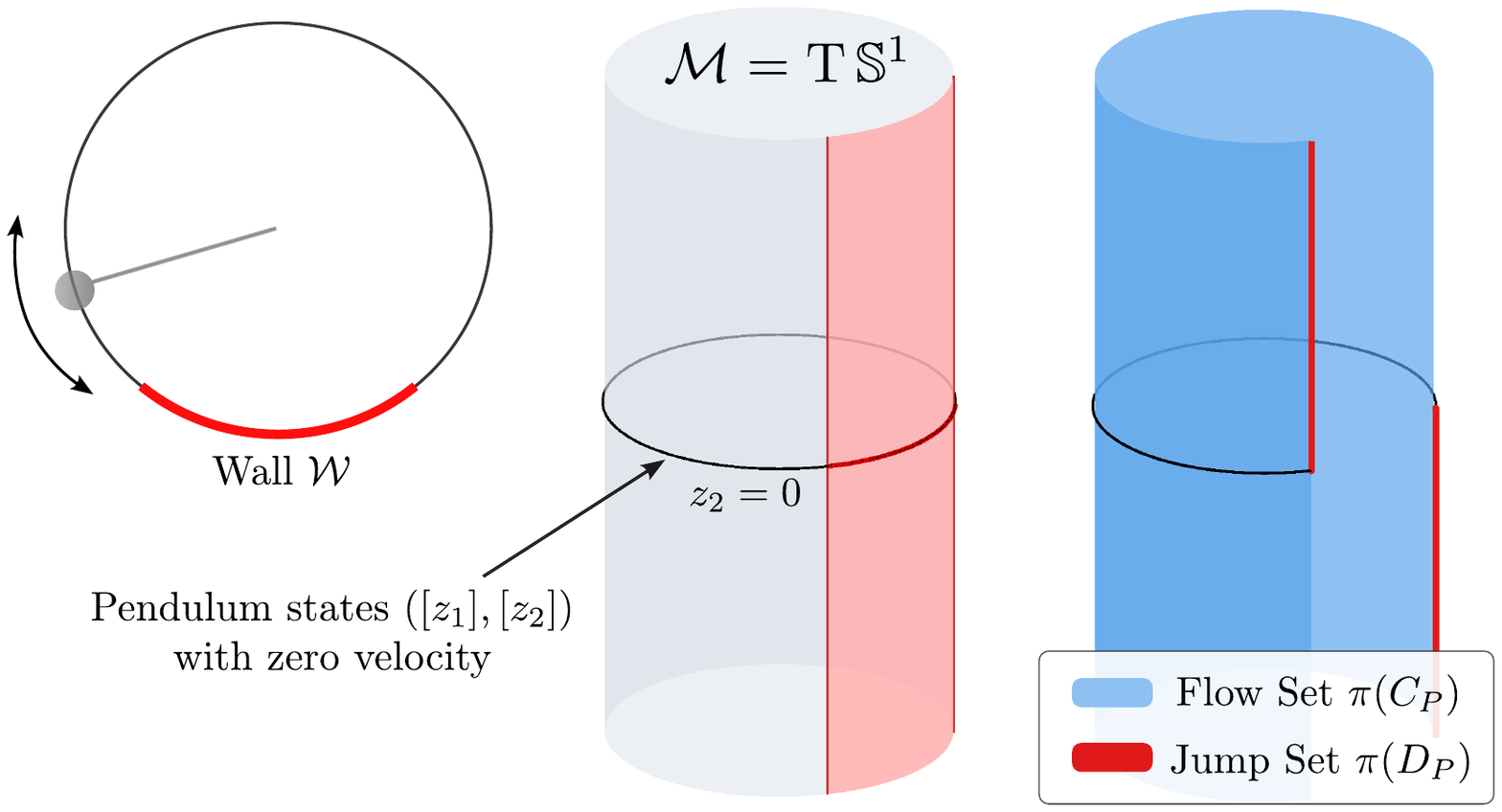}}
        \caption{
            \emph{Left:} Pendulum evolution on the circle $\mathbb{S}^1$, with the wall $\cal W$ depicted in red. 
            \emph{Center:} The state space manifold $\M$. The circle in black denotes the states of the pendulum with zero velocity. The red region denotes those states with arbitrary velocity but position in $\cal W$. 
            \emph{Right:} The projection of $C_P$ and $D_P$ onto $\M$ in blue and red, respectively.
        }
        \label{fig:pendulum-flowAndJump}
    \end{figure}

    On this bundle, the continuous dynamics are
    \begin{align}
        \dot{{z}} = F_P(\chi) \qquad \chi \in C_P,
    \end{align}
    where
    \begin{align*}
        C_P &\coloneqq \left\{ \chi \in \U: {z} = \pi(\chi),
        \pi_{\mathbb{S}^1}({z}) \in \overline{\mathbb{S}^1\setminus \cal W} \right\},
        \\ 
        F_P(\chi) &\coloneqq \mathfrak{S}(\pi(\chi)) + \vlift{\pi(\chi)}{\operatorname{pr}_2(\chi) - b\pi(\chi)} \qquad \forall \chi\in \U,
    \end{align*}
    where $b>0$ is the viscous-friction coefficient, $\operatorname{pr}_2(\chi)$ is the applied force that is the control input, and $-b\pi(\chi)$ is the friction force. The geodesic spray $\mathfrak{S}:\M\to\T{}{\M}$ describes the unforced motion, while the vertical lift $\operatorname{vlft}_x:\T{\pi_{\mathbb{S}^1}(x)}{\mathbb{S}^1}\to\T{x}{\M}$ maps the net force on the pendulum to a tangent vector to $\mathbb{S}^1$. 
    Together, these terms express the second-order pendulum dynamics as first-order dynamics on $\M=\T{}{\mathbb{S}^1}$; see~\cite[Remark~3.103, p.~113]{BulloLewis} and~\cite[Definition~3.73, p.~96]{BulloLewis} for details. Recall from~\eqref{eq:state-using-input} the term $\pi(\chi)$ denotes the state, and remains constant on each fiber $\pi^{-1}(z)$. Therefore, purely state-dependence of terms in the flow set and the flow map is captured by $\pi(\chi)$, whereas the input dependence is captured by the term $\operatorname{pr}_{2}(\chi)$.

    The collision of the pendulum with the wall, which is captured by the forthcoming jump dynamics, occurs at the wall boundary when its velocity points ``into the wall.'' For $\chi=(v,w,\lambda)\in\U$, the third component $\lambda=\operatorname{pr}_3(\chi)\in[0,1]$ is the controlled coefficient of restitution. The jump dynamics are
    \begin{align}
        {z}^+ = G_P(\chi) \qquad \chi \in D_P,
    \end{align}
    where
    \begin{align*}
        D_P &\coloneqq \left\{ \chi \in \U: {z} = \pi(\chi),
        \pi_{\mathbb{S}^1}({z}) \in \partial \cal W, {z} \in \mathrm{T}^{\mathbb{S}^1}_{\cal W}([z_1])\right\},
        \\
        G_P(\chi) &\coloneqq  -\operatorname{pr}_3(\chi) \;\pi(\chi) \qquad \forall \chi\in \U.
    \end{align*}
    The jump set $D_P$ contains the states at the wall boundary whose velocity points into the wall. At such a state, $G_P$ reverses the velocity and scales its magnitude by $\lambda =  \operatorname{pr}_3(\chi) \in [0,1]$, where $\lambda=1$ corresponds to an elastic impact and $\lambda=0$ to a perfectly inelastic impact. Scalar multiplication by $-1$ is well-defined because $\pi(\chi)\in\T{\pi_{\mathbb{S}^1}(\pi(\chi))}{\mathbb{S}^1}$ belongs to the tangent space which is a vector space, and it leaves the position unchanged.

    Combining the above flow and jump dynamics, with ${h_P}({\chi})\coloneqq\pi({\chi})$ for each ${\chi}\in\U$, we obtain the open geometric hybrid system $\hybrid_P=(C_P,F_P,D_P,G_P,{h_P},\mathcal{B})$ as in~\eqref{eq:hybridplant} on the relaxed $C^1$-control bundle $\mathcal{B}=(\M,\U,\pi)$. The above construction of $C_P$ and $D_P$, projected onto $\M$, is illustrated in Figure~\ref{fig:pendulum-flowAndJump}.
\end{example}

\subsubsection{Static State-Feedback Control Laws}

Next, following~\cite{Tabuada2005Quotients}, we define a static state-feedback control law as a section of the control bundle.

\begin{definition}[Static state-feedback control law]
    \label{def:staticFeedback}
    Given a relaxed $C^1$-control bundle $\mathcal{B} = (\M, \U, \pi)$ and an {open geometric hybrid system} $\hybrid_P = (C_P, F_P, D_P, G_P, {h_P}, \mathcal{B})$, the section $\kappa :\M \to \U$ is a \emph{static state-feedback control law}. 
\end{definition}

Note that Definition~\ref{def:staticFeedback} imposes no regularity on $\kappa$. In practice, however, regularity of $\kappa$ is crucial for the existence of solutions to the resulting closed-loop system. In particular, Borel measurability of $\kappa$ is the natural compatibility condition between a static state-feedback law and the solution concept. To this end, the following result, which is a variation of Filippov's lemma~\cite[Lemma~5.2]{goebel_setvaluedAnalysis}, shows that if $\kappa$ is Borel measurable and locally essentially precompact, then it maps a hybrid arc to a hybrid input.

\begin{lemma}
    \label{lem:hybridInput}
    Given a relaxed $C^1$-control bundle $\mathcal{B} = (\M, \U, \pi)$ and an {open geometric hybrid system} $\hybrid_P$, let $\kappa : \M \to \U$ be a static state-feedback control law for $\hybrid_P$. Suppose that $\kappa$ is Borel measurable and locally precompact. Then, for each hybrid arc $z : \dom{z} \to \M$, the composition $\kappa\circ z : \dom{z} \to \U$ is a hybrid input. 
\end{lemma}

\begin{proof}
    Pick any $j\in \mathbb{N}$ such that $I^j_{z}\coloneqq \{t : (t,j)\in \dom{z}\}$ is nonempty. As $z$ is a hybrid arc, $t\mapsto z(t,j)$ is locally absolutely continuous on $I^j_{z}$. Then, it is also continuous and, therefore, Borel measurable; namely, $(\mathfrak{B}(\R{}_{\geq 0}), \mathfrak{B}(\M))$-measurable on $I^j_{z}$. As $\kappa$ is also Borel measurable; namely $(\mathfrak{B}(\M), \mathfrak{B}(\U))$-measurable, and as the composition of Borel measurable functions is also Borel measurable, it follows that $t\mapsto (\kappa\circ z)(t,j)$ is $(\mathfrak{B}(\R{}_{\geq 0}), \mathfrak{B}(\U))$-measurable on $I^{j}_z$. As $\mathfrak{B}(\R{}_{\geq 0}) \subset \mathfrak{L}(\R{}_{\geq 0})$, the previous statement implies that $t\mapsto (\kappa\circ z)(t,j)$ is also $(\mathfrak{L}(\R{}_{\geq 0}),\mathfrak{B}(\U))$-measurable and, therefore, Lebesgue measurable on $I^j_{z}$. Finally, as $\kappa$ is locally precompact and $t\mapsto z(t,j)$ is continuous on $I^{j}_z$, the composition is locally essentially precompact on $I^j_{z}$. As $\dom{z} = \dom(\kappa\circ z)$, $I^{j}_z = I^{j}_{\kappa\circ z}$, and the proof is complete. 
\end{proof}

The local precompactness assumption in Lemma~\ref{lem:hybridInput} is not restrictive as practical applications often require the control inputs to be ``bounded.''

\begin{remark}
    Continuity of a static state-feedback control law $\kappa$ also implies its Borel measurability and local essential precompactness. However, requiring continuity of $\kappa$ may be too strong for several applications; for example, sliding mode control. Borel measurability provides suitable regularity for $\kappa$ while also avoiding stronger-than-necessary assumptions. 
\end{remark}

\section{Barrier Certificates for Set Invariance}
\label{sec:mBF}

The previous section modeled closed and open geometric hybrid systems on manifolds. In this section, we study a closed-loop geometric hybrid system $\hybrid = (C, F, D, G, \M)$ and establish sufficient conditions for forward invariance of a set $S \subset C \cup D$. These conditions generalize~\cite{MAGHENEM2021109328} and are presented using a \emph{multiple barrier function candidate} $B : \dom{B}\to \R{N}$, $N \in \{1,2,\ldots\}$. The function $B$ has $N$ components given by functions $B_i : \dom{B} \to \R{}$, where $i\in \mathcal{I}$ and $\mathcal{I}$ denotes an index set with $\operatorname{card}(\mathcal{I}) = N$. For each $i\in \mathcal{I}$, the zero-sublevel set of $B_i$, restricted to $C\cup D$, defines the set
\begin{align}
    \label{eq:Si_BF}
    S_i \coloneqq \{x\in C\cup D : B_i(x) \leq 0\}. 
\end{align}
The set of interest, namely, $S \subset C\cup D$, is then parametrized by the function $B$ so that $S$ is the intersection of $S_i$'s over $i \in \mathcal{I}$.  

Formally, the above construction of a multiple barrier function candidate is presented below. 

\begin{definition}[Multiple barrier function candidate]
    \label{def:multipleBarrierFunctionCandidate}
    A function $B : \dom{B} \to \R{N}$ is a \emph{multiple barrier function candidate} with respect to a nonempty set $S \subset \M$ for $\hybrid = (C, F, D, G, \M)$ if the following conditions hold:
    \begin{enumerate}
        \item \label{item:mBF-thm-1} $\overline{C} \cup D \cup G(D) \subset \dom{B}$;
        \item \label{item:mBF-thm-2} $S \subset C \cup D$;
        \item \label{item:mBF-thm-3} for each $i\in \mathcal{I}$, there exists an open set $U_i \subset \M$ with $M_i \coloneqq \{x\in \partial S : B_i(x) = 0\} \subset U_i$, such that $B_i$ is locally Lipschitz on an open set containing $(U_i \setminus S_i) \cap \overline{C}$;
        \item \label{item:mBF-thm-4} $S = \displaystyle{\bigcap_{i\in \mathcal{I}}} S_i$, where $S_i$ for each $i\in\mathcal{I}$ is defined in~\eqref{eq:Si_BF}
    \end{enumerate}
    If $N = 1$, we say that $B$ is a \emph{scalar} barrier function candidate with respect to $S$ for $\hybrid$. 
\end{definition}

Note that $B$ is locally Lipschitz in a neighborhood outside $S$ in $\overline{C}$. Consequently, we use (Clarke) generalized derivative from Definition~\ref{def:generalizedDerivative} to compute the rate of change of $B$ during flows in this neighborhood, and at jumps.  

\begin{definition}
    Given a geometric hybrid dynamical system $\hybrid = (C, F, D, G, \M)$ and a function $f : \dom{f} \to \R{}$ with $\overline{C}\cup D\cup G(D) \subset \dom{f}$ that $f$ is locally Lipschitz on a set $U \subset \dom{f}$,
    \begin{enumerate}
        \item the largest rate of change of $f$ along flows in $U$ is given by
        \begin{align*}
            \dot{f}(x) \coloneqq \sup_{v\in F(x)\cap \Tcone{C}{x}} f^{\circ}(x, v) \qquad \forall x\in U \cap C, 
        \end{align*}
        \item the largest value of $f$ after jumps is given by
        \begin{align*}
            f^+(x) \coloneqq \sup_{g\in G(x)} f(g) \qquad \forall x\in D. 
        \end{align*}
    \end{enumerate}
\end{definition}

We present the following main result of this section.
 
\begin{theorem}[Forward pre-invariance using a multiple barrier function]\label{thm:BimpliesForward}
    Given a geometric hybrid system $\mathcal{H} = (C, F, D, G, \cal M)$ as in \eqref{eq:HS} and a nonempty set $S \subset \cal M$, suppose that $S$ is closed and $B : \dom{B}\to\R{N}$ is a multiple barrier function candidate with respect to $S$ for $\cal H$. For each $i\in \cal I$, let $U_i$ be an open neighborhood of $M_i$ as in Definition~\ref{def:multipleBarrierFunctionCandidate}. Then, the set $S$ is forward pre-invariant for $\cal H$ if, for each $i\in \cal I$,
    \begin{subequations}
    \label{eq:BimpliesFI}
      \begin{align}
          \dot{B}_i(x)  & \leq 0 \qquad\qquad  \forall x \in (U_i\backslash S_{i}) \cap C \label{eq:thm1}
          \\
          B^+_i(x)      & \leq 0 \qquad\qquad \forall x \in D \cap S_i  \label{eq:thm2}\\
          G(D \cap S) & \subset C \cup D. \label{eq:thm3}
      \end{align}
      \end{subequations}
\end{theorem}

\begin{proof}
    The proof follows along similar lines to the proof of~\cite[Theorem~1]{MAGHENEM2021109328}. Suppose by contradiction that~\eqref{eq:BimpliesFI} hold but $S$ is not forward pre-invariant for $\cal H$. Consequently, there exists a solution $\phi$ to $\hybrid$ from $S$ that exits the set $S$. Then, one of the following two cases must occur.
    \begin{enumerate}[label = {(\Roman*)}]
        \item \label{item:Case-I}The solution $\phi$ exits $S$ via a jump. Then, there exists $(t,j)\in\dom{\phi}$ with $\phi(t,j)\in S \cap D$ and $\phi(t,j+1) \notin S$.
        \item \label{item:Case-II} The solution $\phi$ exits $S$ via flows. Then, there exist $j \in \mathbb{N}$ and $0 < t_1 < t_2$ such that $[t_1, t_2]\times \{j\} \subset \dom{\phi}$, $\phi(t_1, j) \in \partial S$, $\phi(s, j) \notin S$ for all $s\in (t_1, t_2]$. 
    \end{enumerate}

    Suppose $\phi$ exits $S$ such that~\ref{item:Case-I} holds. Then, let $(t,j)\in\dom{\phi}$ be such that $\phi(t,j)\in S\cap D$ and $\phi(t,j+1)\notin S$. From item~\ref{item:mBF-thm-4} of Definition~\ref{def:multipleBarrierFunctionCandidate}, there exists $i\in \cal I$ such that $\phi(t,j+1)\notin S_i$. Note also from from~\eqref{eq:thm3} that $\phi(t,j+1) \in C\cup D$. Consequently, \eqref{eq:Si_BF} yields that $B_i^+(\phi(t,j)) > 0$, which contradicts~\eqref{eq:thm2}. Therefore, \ref{item:Case-I} does not hold.

    Alternatively, suppose that $\phi$ exits $S$ such that~\ref{item:Case-II} holds. Recall $(t_1, j), (t_2, j)\in \dom{\phi}$ from \ref{item:Case-II}. Due to continuity of $t\mapsto \phi(t,j)$ and item~\ref{item:mBF-thm-4} of Definition~\ref{def:multipleBarrierFunctionCandidate}, we have that, for $t_2$ sufficiently close to $t_1$, there exists $i\in \cal I$ such that 
    \begin{align}
    \label{eq:phiOutsideS}
        \phi(t,j)\in (U_i\setminus S_i) \cap C \qquad \forall t\in (t_1, t_2).
    \end{align}
 
    \noindent
    Using local Lipschitz continuity of $B_i$ and absolute continuity of $\phi$ on $[t_1, t_2]\times\{j\}$, it follows that
    \begin{align}
    \label{eq:Bchange}
        B_i(\phi(t_2, j)) - B_i(\phi(t_1, j)) = \int_{t_1}^{t_2} \frac{d}{dt}B_i(\phi(t, j)) dt.
    \end{align}
    Since $B_i(\phi(t_1,j)) = 0$ as $\phi(t_1,j)\in \partial S$ and $B_i(\phi(t_2, j)) > 0$ as $\phi(t_2, j)\notin S$, the expression in~\eqref{eq:Bchange} is positive. However, using ~\cite[Lemma~1]{JirwankarSanfelice_CDC2025} for the function $B_i$ yields
    \begin{align*}
        \frac{d}{dt}B_i(\phi(t, j)) \leq \dot{B}_i(\phi(t,j)) \quad \text{for almost all } t\in [t_1, t_2].
    \end{align*}
    Since $\dot{B}_i(\phi(t,j))\leq 0$ for each $t\in [t_1, t_2]$ from~\eqref{eq:thm1} and~\eqref{eq:phiOutsideS}, integrating the above equation results in the expression in~\eqref{eq:Bchange} being nonpositive. This is a contradiction. As \ref{item:Case-I} and \ref{item:Case-II} do not hold, $\phi$ must not leave $S$. 
\end{proof}

\begin{remark}
    As the forward invariance property of $S$ depends on the flow vector field immediately outside $S$ (c.f.~\cite[Ch. 4.6]{AubinCellina1984}), the barrier function decrease condition in~\eqref{eq:thm1} is enforced in a neighborhood outside $S_i$. During jumps, however, a solution starting inside $S$ can potentially jump outside $S$, which, therefore, necessitates the condition in~\eqref{eq:thm2} to be enforced on $D \cap S_i$.
\end{remark}

Theorem~\ref{thm:BimpliesForward} can be invoked to certify that sublevel sets of a Lyapunov function with respect to a nonempty closed set $\A\subset \M$ for $\hybrid$ are forward invariant. This result is a generalization of~\cite[Proposition~6.5]{jirwankar2026barrier}. 

\begin{proposition}
    Consider a geometric hybrid system $\cal H = (C, F, D, G, \cal M)$ such that $C\cup D = \cal M$. Given a nonempty, closed set $\cal A \subset \cal M$, suppose there exists a continuous function $V: \cal M \to \R{}_{\geq 0}$ that is locally Lipschitz continuous on an open set containing $\overline{C}$, and that satisfies $V \in \cal PD(\cal A)$. If there exists $s > 0$ such that
    \begin{align}
            \dot{V}(x) \leq 0 && x \in C \cap L_V(s), \label{eq:dotV}\\
            V^+(x) - V(x) \leq 0 && x \in D \cap L_V(s), \label{eq:deltaV}
    \end{align}
    where $L_V(s)\coloneqq \{x\in \M: V(x)\leq s\}$, then, for each $r \in [0, s)$, the $r$-sublevel set $L_V(r)$ is forward pre-invariant for $\cal H$.
\end{proposition}

\begin{proof}
    Pick any $r \in [0, s)$ and let $S \coloneqq L_V(r) = V^{-1}([0,r])$. Since $V$ is continuous and $[0,r]\subset \R{}$ is closed, $S$ is closed. Additionally, as $r < s$, we have $L_V(r)\subset L_V(s)$. Next, let $r' \in (r,s)$ so that $\interior{L_V(r')}$ denotes an open neighborhood of $S$. By construction, $V$ is locally Lipschitz on $(\interior{L_V(r')}\setminus S)\cap \overline{C}$. Then, consider the function $B : \cal M \to \R{}$ such that
    \[
        B(x) \coloneqq V(x) - r \qquad \forall x\in \M
    \]
    and note that $B$ is a scalar barrier function candidate with respect to $S$ for $\hybrid$.
    
    Following~\eqref{eq:dotV} and using $\interior{L_V(r')}\setminus S \subset L_V(s)$, we have that $\dot{B}(x)\leq 0$ for each $x\in (\interior{L_V(r')}\setminus S)\cap C$, and~\eqref{eq:thm1} holds. Following~\eqref{eq:deltaV} and using $S\subset L_V(s)$, we have that $V^+(x) \leq V(x)$ and, therefore, $B^+(x)\leq B(x)$ for each $x\in D\cap S$. Then, for each $x\in D\cap S$, as $B(x) \leq 0$, we have $B^+(x)\leq 0$, and~\eqref{eq:thm2} holds. Finally, as $C\cup D = \cal M$ by assumption,~\eqref{eq:thm3} holds. Therefore, by Theorem~\ref{thm:BimpliesForward}, $L_V(r)$ is forward pre-invariant for $\cal H$.
\end{proof}

Note that Theorem~\ref{thm:BimpliesForward} certifies forward pre-invariance of $S$, but makes no comment about forward invariance, which requires that each maximal solution to $\hybrid$ starting from $K$ is complete. This result is established by using sufficient conditions for existence of solutions to $\hybrid$; namely, the geometric hybrid basic conditions. 

\begin{definition}[{Geometric hybrid basic conditions~\cite[Definition~3]{JirwankarOchoaSanfelice_CDC2025}}]\label{def:HBCs}
    The system $\hybrid = (C, F, D, G, \M)$ is said to satisfy the {\em geometric hybrid basic conditions} if
    \begin{enumerate}[label = (A\arabic*), leftmargin=*]
        \item \label{item:HBC1} $\M$ is a $C^1$-manifold;
        \item \label{item:HBC2} $C $ and $D$ are subsets of $\M$ that are closed relative to~$\M$;
        \item \label{HBC:F}$F:\M \rightrightarrows \T{}{\M}$ is outer semicontinuous and locally~precompact relative to $C$, $F(x)$ is convex for all $x\in C$, and $C \subset \dom{F}$; 
        \item \label{HBC:G} $G:\M \rightrightarrows \M$ is outer semicontinuous and locally~precompact relative to $D$, and $D \subset \dom{G}$.
    \end{enumerate}
\end{definition}

Under geometric hybrid basic conditions and the assumptions of Theorem~\ref{thm:BimpliesForward}, we can leverage~\cite[Proposition~2]{JirwankarSanfelice_CDC2025} to present conditions on the set $S$ and the data of $\hybrid$ that strengthen the certificate in Theorem~\ref{thm:BimpliesForward} to forward invariance. This result is from~\cite[Corollary~6.7]{jirwankar2026barrier}.

\begin{corollary}
    \label{cor:forwardInvariance}
    Under the assumptions of Theorem~\ref{thm:BimpliesForward}, suppose that $\cal H = (C, F, D, G, \cal M)$ satisfies the geometric hybrid basic conditions in Definition~\ref{def:HBCs}. Then, a nonempty closed set $S\subset \cal M$ is forward invariant for $\cal H$ if
    \begin{enumerate}
        \item the conditions in \eqref{eq:thm1}-\eqref{eq:thm3} hold,
        \item\label{item:VC} for each $\xi\in (S \cap \partial C)\setminus D$, there exists an open neighborhood $U$ of $\xi$ such that
        \[
            F(x) \cap \Tcone{C}{x} \neq \varnothing \qquad \forall x\in U\cap C,
        \]
        \item\label{item:noEscapeSolutions} at least one of the following conditions holds:
        \begin{enumerate}
            \item $S \cap C$ is compact;
            \item The constrained differential inclusion $\dot{x} \in F(x)\;\; x\in C$ does not have a maximal solution with a finite escape time.
        \end{enumerate}
    \end{enumerate}
\end{corollary}

\section{Multiple Control Barrier Functions}

With the goal of designing control laws that render a set forward invariant, we now develop multiple control barrier functions (mCBFs) for the open geometric hybrid system $\hybrid_P=(C_P,F_P,D_P,G_P,{h_P},\mathcal{B})$, where $\mathcal{B}=(\M,\U,\pi)$ is a relaxed $C^1$-control. The construction extends the nonsmooth pointwise conditions in~\cite[Def.~3.1]{montenegro2025pointwise} to the geometric setting of control bundles. Similar to the approach in Section~\ref{sec:mBF}, we begin by parametrizing the set of interest, and then introduce the mCBF conditions and their admissible feedback maps, and finally use Theorem~\ref{thm:BimpliesForward} to establish forward pre-invariance of the resulting closed-loop system.

Let $K\subset\pi(C_P)\cup\pi(D_P)$ be nonempty, let $N\in\{1,2,\ldots\}$, and let $\mathcal{I}$ be a finite index set satisfying $\operatorname{card}(\mathcal{I})=N$. Consider a function $B:\dom{B}\to\R{N}$ with component functions $B_i:\dom{B}\to\R{}$ for each $i\in\mathcal{I}$. Following the multiple-barrier construction in~\cite[Definition~6.1]{jirwankar2026barrier}, define
\begin{align}
    \label{eq:Ki}
    K_i \coloneqq \{x \in \pi(C_P) \cup \pi(D_P) : B_i(x)\leq 0\}
\end{align}
for each $i\in\mathcal{I}$. The mCBF-candidate definition below requires $K$ to be the intersection of these componentwise zero-sublevel sets.

\begin{definition}[mCBF candidate]
    \label{def:mCBF_candidate}
    Given a relaxed $C^1$-control bundle $\mathcal{B} = (\M, \U, \pi)$, an open geometric hybrid system $\hybrid_P = (C_P, F_P, D_P, G_P, {h_P}, \mathcal{B})$, and a nonempty set $K \subset \M$, the function $B : \dom{B} \to \R{N}$ is a \emph{multiple control barrier function (mCBF) candidate for $\hybrid_P$ with respect to $K$} if the following conditions hold:
    \begin{enumerate}
        \item \label{item:1-mCBFcandidate} $\overline{\pi(C_P)} \cup \pi(D_P) \cup G_P(D_P) \subset \dom{B}$;
        \item $K \subset \pi(C_P) \cup \pi(D_P)$ and $K$ is closed;
        \item \label{item:2-mCBFcandidate} for each $i\in \mathcal{I}$, $B_i$ is locally Lipschitz on an open neighborhood of $\overline{\pi(C_P)}$;
        \item \label{item:3-mCBFcandidate} $K = \displaystyle \bigcap_{i\in\mathcal{I}} K_i$, where $K_i$ is defined using $B_i$ in~\eqref{eq:Ki}.
    \end{enumerate}
\end{definition}

\begin{example}[Pendulum with impacts, revisited]
    \label{ex:pendulum-mCBF-candidate}
    Recall the pendulum dynamics in Example~\ref{ex:pendulum-plant-dynamics}. Suppose that the gravitational potential $U : \mathbb{S}^1 \to \R{}_{\geq 0}$ for the pendulum is defined as
    \(
        U([z_1]) \coloneqq \frac{a}{2\pi}(1 - \cos(2\pi z_1))\)
    for each $z_1 \in [0,1]$, where $a > 0$. The total energy of the pendulum is given by $E : \cal M \to \R{}_{\geq 0}$ such that
    \[
        E(x) \coloneqq U([z_1]) + \frac{1}{2} g_{\mathbb{S}^1}^{[z_1]}(x,x) \qquad \forall x = ([z_1], [z_2]) \in \M = \T{}{\mathbb{S}^1},
    \]
    where $[z_1] = \pi_{\mathbb{S}^1}(x)$ denotes the pendulum position on the circle, and $g_{\mathbb{S}^1}^{[z_1]}$ denotes the evaluation of the Riemannian metric $g_{\mathbb{S}^1}$ at $[z_1]$.
    Note that, when the pendulum position is at the boundary of the wall $\W$ with zero velocity, $x\in \partial \cal W \times \{0\} \subset \M$. Due to the definition of $\cal W$ and $U$, it follows that $E(\partial \cal W \times \{0\}) \subset \R{}_{\geq 0}$ is a singleton. Let $c_0 \coloneqq E(\partial \cal W \times \{0\})$, and pick $c_{-1}, c_{1} > 0$ satisfying $c_0 < c_{-1} < c_{1}$. 

    Define $\mathcal{I} \coloneqq \{-1,1\}$. For each $i\in \mathcal{I}$, define $B_i : \M \to \R{}$ as
    \[
        B_i(x) \coloneqq i (E(x) - c_{i}) \qquad \forall x \in \M
    \]
    that, following~\eqref{eq:Ki}, defines the set $K_i$. For each $i\in\mathcal{I}$, as $B_i$ is continuous, the set $\widetilde{K}_i \coloneqq \{x\in \M : B_i(x)\leq 0\}$ is closed. Additionally, from the definitions of $C_P,D_P\subset\U$ in Example~\ref{ex:pendulum-plant-dynamics}, $\pi(C_P) \cup \pi(D_P) = \{x\in \M : \pi_{\mathbb{S}^1}(x) \in \overline{\mathbb{S}^1\setminus \cal W}\}$ is closed. Then, the set
    \[ 
        K = \{x\in \pi(C_P) \cup \pi(D_P) : c_{-1} \leq E(x) \leq c_1\}
    \]
    is closed as well since $K = (\pi(C_P)\cup \pi(D_P))\cap L$, where $L \coloneqq \{x\in \M : E(x)\in [c_{-1}, c_{1}]\}$ is closed as $E$ is continuous. Next, we define $B : \M \to \R{2}$ by $B(x) \coloneqq \begin{pmatrix}
        B_{-1}(x), & B_{1}(x)
    \end{pmatrix}$ for each $x\in \M$. As $\dom{B} = \M$, $B$ is smooth, and $K = \cap_{i\in\mathcal{I}} K_i$ is closed, $B$ is a mCBF candidate for $\hybrid_P$ with respect to $K$. 
\end{example}

To formulate the mCBF inequalities, we follow~\cite{montenegro2025pointwise,Isaly2024Feasibility} and introduce performance functions $\gamma=(\gamma_C,\gamma_D)$, where $\gamma_C:\pi(C_P)\to\R{N}$ and $\gamma_D:\pi(D_P)\to\R{N}$. Their required signs are imposed only on the sets where the flow and jump inequalities are used to certify forward pre-invariance.

\begin{assumption}
    \label{ass:perfomanceFunctions}
    Let $\hybrid_P = (C_P, F_P, D_P, G_P, {h_P}, \mathcal{B})$ be an open geometric hybrid system, with $\mathcal{B} = (\M, \U, \pi)$ a relaxed $C^1$-control bundle, and $B : \dom{B} \to \R{N}$ be a mCBF candidate for $\hybrid_P$ with respect to a nonempty set $K \subset \M$. Then, there exist functions $\gamma_C : \pi(C_P) \to \R{N}$ and $\gamma_{D} : \pi(D_P) \to \R{N}$ and, for each $i\in \mathcal{I}$, an open set $U_i \subset \M$ containing $M_i \coloneqq \{x\in\partial K : B_i(x) = 0\}$, such that
    \begin{align}
        \label{eq:perfomanceFunctions_positive}
        \begin{array}{rl}
            \gamma_C^{(i)}(x) \geq 0 & \qquad  \forall x\in ({U}_i\setminus K_i) \cap \pi(C_P),
            \\[5pt]
            \gamma_D^{(i)}(x) \geq 0 & \qquad \forall x\in K_i \cap \pi(D_P),
        \end{array}
    \end{align}
    holds for each $i\in\mathcal{I}$, where $\gamma_C^{(i)}$ and $\gamma_D^{(i)}$ denote the $i$'th component of $\gamma_C$ and $\gamma_D$, respectively. 
\end{assumption}

\begin{definition}[mCBF]
    \label{def:mCBF}
    Given a relaxed $C^1$-control bundle $\mathcal{B} = (\M, \U, \pi)$, an open geometric hybrid system $\hybrid_P = (C_P, F_P, D_P, G_P, {h_P}, \mathcal{B})$, a nonempty set $K$, a mCBF candidate $B : \dom{B}\to \R{N}$ for $\hybrid_P$ with respect to $K$, and performance functions $\gamma_C : \pi(C_P) \to \R{N}$ and $\gamma_D : \pi(D_P) \to \R{N}$ that satisfy Assumption~\ref{ass:perfomanceFunctions}, the function $B$ is a \emph{multiple control barrier function (mCBF) for $\hybrid_P$ with respect to $(K, \gamma)$} if the following conditions hold:
    \begin{enumerate}
        \item \label{item:1-mCBF} $ C_P \subset \dom{F_P}$ and $D_P \subset \dom{G_P}$;
        \item \label{item:2-mCBF} for each $x\in \pi(C_P)$, there exists $\chi\in \pi^{-1}(x)\cap C_P$ such that
        \[
            \sup_{f \in F_P(\chi)} B_i^{\circ}(x, f) \leq -\gamma_C^{(i)}(x) \qquad \forall i\in \mathcal{I};
        \]
        \item \label{item:3-mCBF} for each $x\in \pi(D_P)$, there exists $\chi\in \pi^{-1}(x) \cap D_P$ such that
        \[
            \sup_{g \in G_P(\chi)} B_i(g) \leq -\gamma_D^{(i)}(x) \qquad \forall i\in\mathcal{I}.
        \]
    \end{enumerate}
\end{definition}

\begin{remark}
    \label{remark:fullDomains}
    Definition~\ref{def:mCBF} requires the pointwise mCBF constraints to admit at least one input at each state during flows and jumps. This requirement does not impose a decrease condition on the mCBF inside $K$ during flows. Indeed, Assumption~\ref{ass:perfomanceFunctions} does not require $\gamma_C^{(i)}$ to be nonnegative at points in $K$ that do not belong to $(U_i\setminus K_i)\cap\pi(C_P)$ for any $i\in\mathcal{I}$. At each such point, $\gamma_C^{(i)}$ may be assigned a sufficiently large negative value. Consequently, $-\gamma_C^{(i)}$ is positive, and the flow constraint permits the components of $B$ to increase while remaining bounded above by $-\gamma_C$. Such design of $\gamma_C^{(i)}$ can make $\mathbb{U}_C(x)$ nonempty for each $x\in \pi(C_P)$ without restricting the behavior of the mCBF inside $K$.
\end{remark}

Given a mCBF $B$ for $\hybrid_P$ with respect to $(K,\gamma)$, define $\mathbb{U}_C : \pi(C_P) \setto \U$ and $\mathbb{U}_D : \pi(D_P) \setto \U$ to collect the flow and jump inputs, respectively, that satisfy all componentwise mCBF inequalities. In particular, let
\begin{align}\label{eq:bbU}
    \begin{array}{lr}
        \mathbb{U}_C(x) \coloneqq \left \{\chi \in \pi^{-1}(x) \cap C_P :  \underset{f\in F_P(\chi)}{\sup} B_{i}^{\circ}(x, f) \leq -\gamma_C^{(i)}(x) \quad \forall i\in \mathcal{I} \right \} & \forall x \in \pi(C_P),
        \\[10pt]
        \mathbb{U}_D(x) \coloneqq \left\{\chi \in \pi^{-1}(x) \cap D_P : \underset{g\in G_P(\chi)}{\sup} B_{i}(g) \leq -\gamma_D^{(i)}(x) \quad \forall i\in \mathcal{I}\right\} & \forall x \in \pi(D_P). 
    \end{array}
\end{align}

The following result establishes that $\mathbb{U}_C$ and $\mathbb{U}_D$ have nonempty values on $\pi(C_P)$ and $\pi(D_P)$, respectively. 

\begin{lemma}
    \label{lem:mCBFFullDomains}
    Given a relaxed $C^1$-control bundle $\mathcal{B} = (\M, \U, \pi)$, an open geometric hybrid system $\hybrid_P = (C_P, F_P, D_P, G_P, {h_P}, \mathcal{B})$, a nonempty set $K$, a mCBF candidate $B : \dom{B}\to \R{N}$ for $\hybrid_P$ with respect to $K$, and performance functions $\gamma_C : \pi(C_P) \to \R{N}$ and $\gamma_D : \pi(D_P) \to \R{N}$ that satisfy Assumption~\ref{ass:perfomanceFunctions}, suppose that $B$ is a \emph{multiple control barrier function (mCBF) for $\hybrid_P$ with respect to $(K, \gamma)$}. Then, $\dom{\mathbb{U}_C} = \pi(C_P)$ and $\dom{\mathbb{U}_D} = \pi(D_P)$. 
\end{lemma}

\begin{proof}
    The conclusion follows immediately from Definition~\ref{def:mCBF} and using the fact that $\mathbb{U}_C(x)$ (resp., $\mathbb{U}_D(x)$) is nonempty for each $x\in \pi(C_P)$ (resp., $\pi(D_P)$). 
\end{proof}

Applying the set-valued maps $\mathbb{U}_C$ and $\mathbb{U}_D$ as feedback laws to $\hybrid_P$ yields the closed-loop geometric hybrid system $\hybrid_{\mathbb{U}}=(C_{\mathbb{U}},F_{\mathbb{U}},D_{\mathbb{U}},G_{\mathbb{U}},\M)$ given by
\begin{align}
    \label{eq:closedLoopU}
    \hybrid_{\mathbb{U}} : \left\{
        \begin{array}{ll}
            \dot{x} \in F_{\mathbb{U}}(x)
            \coloneqq \displaystyle\bigcup_{\chi\in\mathbb{U}_C(x)} F_P(\chi)
            & \quad x\in C_{\mathbb{U}}\coloneqq\dom{\mathbb{U}_C},
            \\[8pt]
            x^+ \in G_{\mathbb{U}}(x)
            \coloneqq \displaystyle\bigcup_{\chi\in\mathbb{U}_D(x)} G_P(\chi)
            & \quad x\in D_{\mathbb{U}}\coloneqq\dom{\mathbb{U}_D}.
        \end{array}
    \right.
\end{align}

The following result proves forward pre-invariance of $K$ for $\hybrid_{\mathbb{U}}$. When a control law is chosen as a suitable selection from $\mathbb{U}_C$ and $\mathbb{U}_D$, additional conditions for forward pre-invariance of $K$ are presented in forthcoming Lemma~\ref{lem:singleSelectionInvariance}. 

\begin{theorem}[Forward pre-invariance under mCBF feedbacks]
    \label{thm:mCBFForwardPreInvariant}
    Given a relaxed $C^1$-control bundle $\mathcal{B}=(\M,\U,\pi)$, an open geometric hybrid system $\hybrid_P=(C_P,F_P,D_P,G_P,{h_P},\mathcal{B})$, a nonempty set $K\subset\M$, and performance functions $\gamma=(\gamma_C,\gamma_D)$ satisfying Assumption~\ref{ass:perfomanceFunctions}, let $B:\dom{B}\to\R{N}$ be a mCBF for $\hybrid_P$ with respect to $(K,\gamma)$. 
    Suppose that \begin{align}
        \label{eqref:GU_in_CUDU}
        G_{\mathbb{U}}(K\cap D_{\mathbb{U}})\subset C_{\mathbb{U}}\cup D_{\mathbb{U}}.
    \end{align}
    Then, $K$ is forward pre-invariant for the closed-loop system $\hybrid_{\mathbb{U}}$ in~\eqref{eq:closedLoopU}.
\end{theorem}

\begin{proof}
    For each $i\in\mathcal{I}$, define $S_i\coloneqq\{x\in C_{\mathbb{U}}\cup D_{\mathbb{U}}:B_i(x)\leq0\}$.
    We verify the hypotheses of Theorem~\ref{thm:BimpliesForward} for $\hybrid_{\mathbb{U}}$, with the set $S$ therein denoted by $K$.
    By Lemma~\ref{lem:mCBFFullDomains},
    $C_{\mathbb{U}}=\pi(C_P)$ and $D_{\mathbb{U}}=\pi(D_P)$.
    Therefore, item~\ref{item:1-mCBFcandidate} of Definition~\ref{def:mCBF_candidate} yields
    \begin{align*}
        \overline{C_{\mathbb{U}}}\cup D_{\mathbb{U}}\cup G_{\mathbb{U}}(D_{\mathbb{U}})
        \subset
        \overline{\pi(C_P)}\cup\pi(D_P)\cup G_P(D_P)
        \subset\dom{B},
    \end{align*}
    where $G_{\mathbb{U}}(D_{\mathbb{U}})\subset G_P(D_P)$ follows because, for each $x\in D_{\mathbb{U}}$ and each $\chi\in\mathbb{U}_D(x)$, $G_P(\chi)\subset G_P(D_P)$. Lemma~\ref{lem:mCBFFullDomains} and Definition~\ref{def:mCBF_candidate} yield $K\subset C_{\mathbb{U}}\cup D_{\mathbb{U}}$.
    Additionally, for each $i\in\mathcal{I}$, $S_i=K_i$, and $B_i$ is locally Lipschitz on $(U_i\setminus S_i)\cap\overline{C_{\mathbb{U}}} =(U_i\setminus K_i)\cap\overline{\pi(C_P)}$
    by item~\ref{item:2-mCBFcandidate} of Definition~\ref{def:mCBF_candidate}. Next, item~\ref{item:3-mCBFcandidate} of Definition~\ref{def:mCBF_candidate} yields $K=\bigcap_{i\in\mathcal{I}}S_i.$
    Thus, $B$ is a multiple barrier function candidate for $\hybrid_{\mathbb{U}}$ with respect to $K$.

    Pick any $i\in\mathcal{I}$ and any $x\in(U_i\setminus S_i)\cap C_{\mathbb{U}}$. Since $S_i=K_i$, for each $\chi\in\mathbb{U}_C(x)$, the definition of $\mathbb{U}_C$ and Assumption~\ref{ass:perfomanceFunctions} yield
    \begin{align*}
        \sup_{f\in F_P(\chi)}B_i^{\circ}(x,f)
        \leq-\gamma_C^{(i)}(x)
        \leq 0.
    \end{align*}
    Since the above inequality holds for each $\chi\in\mathbb{U}_C(x)$, it follows that
    \[
        \sup_{f\in \bigcup_{\chi\in\mathbb{U}_C(x)}F_P(\chi)}B_i^{\circ}(x,f) = \sup_{f\in {F}_{\mathbb{U}}(x)} B_{i}^{\circ}(x,f) \leq 0. 
    \]

    Next, pick any $i\in\mathcal{I}$ and any $x\in S_i\cap D_{\mathbb{U}}$. Since $S_i=K_i$, for each $\chi\in\mathbb{U}_D(x)$, the definition of $\mathbb{U}_D$ and Assumption~\ref{ass:perfomanceFunctions} yield
    \(
        \sup_{g\in G_P(\chi)}B_i(g)
        \leq-\gamma_D^{(i)}(x)
        \leq0.
    \)
    Consequently, $\sup_{g\in G_{\mathbb{U}}(x)}B_i(g)\leq0$.
    Finally, using~\eqref{eqref:GU_in_CUDU}, Theorem~\ref{thm:BimpliesForward} implies that $K$ is forward pre-invariant for $\hybrid_{\mathbb{U}}$.
\end{proof}

Theorem~\ref{thm:mCBFForwardPreInvariant} proves forward pre-invariance of $K$ using the set-valued feedback $\mathbb{U}_C$ and $\mathbb{U}_D$. However, it makes no comment on the existence of nontrivial solutions to the closed-loop system $\hybrid_{\mathbb{U}}$. Indeed, in the absence of any solutions to $\hybrid_{\mathbb{U}}$ from $K$, forward pre-invariance of $K$ holds vacuously. To strengthen this result, we present the following sufficient conditions that ensure existence of nontrivial solutions to $\hybrid_{\mathbb{U}}$; see~\cite[Proposition~1]{JirwankarOchoaSanfelice_CDC2025}.

\begin{corollary}[Existence of solutions under mCBF feedback]
    \label{cor:mCBFExistenceCompleteness}
    Under the assumptions of Theorem~\ref{thm:mCBFForwardPreInvariant}, suppose that $\hybrid_{\mathbb U}$ satisfies the geometric hybrid basic conditions. If, additionally, for each $\xi\in (K\cap\partial C_{\mathbb U})\setminus D_{\mathbb U}$, there exists an open neighborhood $V_\xi$ of $\xi$ such that
    \begin{align}
       F_{\mathbb U}(x)\cap\Tcone{C_{\mathbb U}}{x}\neq\varnothing
       \qquad \forall x\in V_\xi\cap C_{\mathbb U},
       \label{eq:mCBFViability}
    \end{align}
    then, for each $\xi\in K$, there exists a nontrivial solution to $\hybrid_{\mathbb U}$ from $\xi$. If, additionally, either $K\cap C_{\mathbb U}$ is compact or the constrained differential inclusion $\dot x\in F_{\mathbb U}(x)$  $x\in C_{\mathbb U}$ has no maximal solution with finite escape time, then each maximal solution from $K$ is complete and $K$ is forward invariant for $\hybrid_{\mathbb U}$.
\end{corollary}

\begin{proof}
    Theorem~\ref{thm:mCBFForwardPreInvariant} yields that $K$ is forward pre-invariant. Then, using~\eqref{eq:mCBFViability} with~\cite[Proposition~1]{JirwankarSanfelice_CDC2025} and noting that $\Tcone{C_{\mathbb{U}}}{x} = \T{x}{\M}$ for each $x\in K\cap \interior{C_{\mathbb{U}}}$, the geometric hybrid basic conditions yield a nontrivial solution to $\hybrid_{\mathbb{U}}$ from each $\xi\in K$. Finally, if $K \cap C_{\mathbb{U}}$ is compact or the constrained differential inclusion has no solutions that exhibit finite escape times, then items b and c of~\cite[Proposition~1]{JirwankarSanfelice_CDC2025} do not hold, causing each maximal solution to $\hybrid_{\mathbb{U}}$ from $K$ to be complete. This proves forward invariance of $K$ for $\hybrid_{\mathbb{U}}$. 
\end{proof}

Theorem~\ref{thm:mCBFForwardPreInvariant} permits the set-valued flow and jump feedbacks to select different inputs at a state in $\pi(C_P)\cap\pi(D_P)$ so that the pointwise conditions on the mCBF hold with the chosen input. A single static state-feedback control law $\kappa : \pi(C_P) \cup \pi(D_P) \to \U$, however, assigns a unique input to each state. Consequently, at a state where both regimes are possible, that input must satisfy both the flow and jump mCBF constraints. In particular, given the set-valued maps $\mathbb{U}_C$ and $\mathbb{U}_D$ as in~\eqref{eq:bbU}, define the set-valued map
\[
    \Gamma : \pi(C_P) \cup \pi(D_P) \setto \U
\] 
as
\begin{align}
    \label{eq:Gamma}
    \Gamma(x) \coloneqq \left\{
        \begin{array}{cl}
            \mathbb{U}_C(x) & \quad \text{if } x \in \pi(C_P)\setminus \pi(D_P)
            \\
            \mathbb{U}_D(x) & \quad \text{if } x \in \pi(D_P)\setminus \pi(C_P)
            \\
            \mathbb{U}_C(x) \cap \mathbb{U}_D(x) & \quad \text{if } x \in \pi(C_P)\cap \pi(D_P)
        \end{array}
    \right. \qquad \forall x \in \pi(C_P) \cup \pi(D_P). 
\end{align}
Then, we choose the control law $\kappa : \pi(C_P) \cup \pi(D_P) \to \U$ to be a selection from $\Gamma$, namely, $\kappa(x)\in \Gamma(x)$ for each $x\in \pi(C_P)\cup \pi(D_{P})$. For such a selection, the following result establishes forward pre-invariance of $K$ for the resulting closed-loop system $\hybrid_{\kappa} = (C_{\kappa}, F_{\kappa}, D_{\kappa}, G_{\kappa}, \M)$, where $C_{\kappa} \coloneqq \{x \in \pi(C_P) \cup \pi(D_P) : \kappa(x)\in C_P\}$, $D_{\kappa} \coloneqq \{x\in \pi(D_P) : \kappa(x) \in D_P\}$, $F_{\kappa}(x) \coloneqq F_P(\kappa(x))$ for each $x\in \dom{\kappa}$, and $G_{\kappa}(x) \coloneqq G_P(\kappa(x))$ for each $x\in\dom{\kappa}$. Following the notation in~\cite{montenegro2025pointwise}, we say that such a selection $\kappa$ is \emph{safeguarding}. 

\begin{lemma}[Forward pre-invariance under a single selection]
    \label{lem:singleSelectionInvariance}
    Under the assumptions of Theorem~\ref{thm:mCBFForwardPreInvariant} and with $\Gamma$ as defined in~\eqref{eq:Gamma}, suppose that $\dom{\Gamma} = \pi(C_P) \cup \pi(D_P)$. Then, for each selection $\kappa:\pi(C_P) \cup \pi(D_P) \to \U$ from $\Gamma$, the resulting closed-loop system $\hybrid_{\kappa} = (C_{\kappa}, F_{\kappa}, D_{\kappa}, G_{\kappa}, \M)$ is such that $C_{\kappa} = \pi(C_P)$, $D_{\kappa} = \pi(D_P)$, and $K$ is forward pre-invariant for $\hybrid_{\kappa}$. 
\end{lemma}

\begin{proof}
    As $\dom{\Gamma} = \pi(C_P) \cup \pi(D_P)$, $\Gamma(x)\neq \varnothing$ for each $x\in \pi(C_P) \cup \pi(D_P)$. Pick any selection $\kappa : \dom{\Gamma} \to \U$ from $\Gamma$. Then, 
    pick $x\in\pi(C_P)$, which causes $\kappa(x) \in \Gamma(x)$. If $x\notin \pi(D_P)$, then $\kappa(x)\in \mathbb{U}_C(x)$, and if $x\in\pi(D_P)$, then $\kappa(x)\in \mathbb{U}_C(x) \cap \mathbb{U}_D(x) \subset \mathbb{U}_C(x)$. Consequently, for each $x\in \pi(C_P)$, $\kappa(x) \in \mathbb{U}_C(x) \subset \pi^{-1}(x)\cap C_P \subset C_P$. As this holds for each $x\in \pi(C_P)$, we have that $\pi(C_P) \subset C_{\kappa}$. Conversely, for each $x\in C_{\kappa}$, $\kappa(x)\in C_P$. Since $\kappa$ is a section
    of $\pi$, $\pi(\kappa(x))=x$, which yields $x\in\pi(C_P)$. Therefore, $C_\kappa=\pi(C_P)=C_{\mathbb U}$.
    A similar reasoning yields $D_{\kappa} = \pi(D_P) = D_{\mathbb{U}}$. 

    Next, for each $x\in C_\kappa$,~\eqref{eq:Gamma} and the
    definition of $F_{\mathbb U}$ yield $F_\kappa(x)=F_P(\kappa(x))\subset F_{\mathbb U}(x)$. 
    Similarly, for each $x\in D_\kappa$, we have that $G_\kappa(x)=G_P(\kappa(x))\subset G_{\mathbb U}(x)$. 
    Therefore, each solution to $\hybrid_\kappa$ is a solution to $\hybrid_{\mathbb U}$. Theorem~\ref{thm:mCBFForwardPreInvariant} establishes that each solution to $\hybrid_{\mathbb U}$ starting from
    $K$ remains in $K$. Consequently, each solution to $\hybrid_\kappa$ starting from $K$ remains in $K$, and $K$ is forward pre-invariant for $\hybrid_\kappa$.
\end{proof}

\begin{remark}
    If the selection $\kappa$ from $\Gamma$ does not have sufficient regularity, the conclusion of Lemma~\ref{lem:singleSelectionInvariance} holds vacuously. Stronger regularity conditions on the data of the closed-loop system are required to avoid vacuous statements. To this end, the next section identifies conditions under which $\Gamma$ admits a Borel-measurable pointwise optimal selection.
\end{remark}

Corollary~\ref{cor:mCBFExistenceCompleteness} and Lemma~\ref{lem:singleSelectionInvariance} require the assumptions of Theorem~\ref{thm:mCBFForwardPreInvariant}. However, checking whether~\eqref{eqref:GU_in_CUDU} holds is a nontrivial task as it involves the functions $\mathbb{U}_C$ and $\mathbb{U}_D$, whose expressions may be difficult to compute. Instead, we present the following sufficient condition on the data of the open geometric hybrid system $\hybrid_P$ that ensures satisfaction of~\eqref{eqref:GU_in_CUDU}.

\begin{lemma}
    \label{lem:plantPostJumpCondition}
    Given a relaxed $C^1$-control bundle $\mathcal{B}=(\M,\U,\pi)$, an open geometric hybrid system $\hybrid_P=(C_P,F_P,D_P,G_P,{h_P},\mathcal{B})$, a nonempty set $K\subset\M$, and performance functions $\gamma_C : \pi(C_P) \to \R{N}$ and $\gamma_D : \pi(D_P) \to \R{N}$ that satisfy Assumption~\ref{ass:perfomanceFunctions}, let $B : \dom{B}\to \R{N}$ denote a mCBF for $\hybrid_P$ with respect to $(K, \gamma)$. Suppose that $G_P(D_P)\subset \pi(C_P)\cup\pi(D_P)$. Then, the closed-loop system $\hybrid_{\mathbb{U}}$ in~\eqref{eq:closedLoopU} satisfies $G_{\mathbb{U}}(K\cap D_{\mathbb{U}})
        \subset C_{\mathbb{U}}\cup D_{\mathbb{U}}$.
\end{lemma}

\begin{proof}
    Let $x\in K\cap D_{\mathbb{U}}$ and $g\in G_{\mathbb{U}}(x)$. By the definition of $G_{\mathbb{U}}$, there exists $\chi\in\mathbb{U}_D(x)$ such that $g\in G_P(\chi)$. The definition of $\mathbb{U}_D$ causes $\chi
    \in D_P$ and, therefore,
    \begin{align*}
        g\in G_P(\chi)\subset G_P(D_P)
        \subset \pi(C_P)\cup\pi(D_P),
    \end{align*}
    where the last inclusion follows by assumption. Additionally, Lemma~\ref{lem:mCBFFullDomains} and the definitions in~\eqref{eq:closedLoopU} yield $C_{\mathbb{U}}=\pi(C_P)$ and $D_{\mathbb{U}}=\pi(D_P)$.
    Therefore, $g\in C_{\mathbb{U}}\cup D_{\mathbb{U}}$. Since $x$ and $g$ were arbitrary, the proof is complete.
\end{proof}

Following Corollary~\ref{cor:mCBFExistenceCompleteness}, a nonvacuous invariance certificate is obtained by asserting geometric hybrid basic conditions for the closed-loop system $\hybrid_{\mathbb{U}}$. As these conditions rely on the data of the open geometric hybrid system $\hybrid_P$ and the set-valued control laws $\mathbb{U}_C$ and $\mathbb{U}_D$, the following conditions ensure satisfaction of said geometric hybrid basic conditions.

\begin{proposition}
    \label{prop:GHBCs}
    Given a relaxed $C^1$-control bundle $\mathcal{B}=(\M,\U,\pi)$, an open geometric hybrid system $\hybrid_P=(C_P,F_P,D_P,G_P,{h_P},\mathcal{B})$, a nonempty closed set $K \subset \M$, and performance functions $\gamma=(\gamma_C,\gamma_D)$ satisfying Assumption~\ref{ass:perfomanceFunctions}, let $B:\dom{B}\to\R{N}$ be a mCBF for $\hybrid_P$ with respect to $(K,\gamma)$. Then, with the set-valued maps $\mathbb{U}_C : \dom{\mathbb{U}_C} \setto \U$ and $\mathbb{U}_D : \dom{\mathbb{U}_D}\setto \U$ in~\eqref{eq:bbU}, suppose that the following conditions hold:
    \begin{enumerate}
        \item \label{item:1-prop} $\pi(C_P)$ and $\pi(D_P)$ are closed;
        \item \label{item:2-prop} $F_P$ is outer semicontinuous and locally precompact relative to $C_P$, $F_P(\chi)$ is convex for each $\chi\in C_P$, and $C_P \subset \dom{F_P}$;
        \item \label{item:3-prop} $G_P$ is outer semicontinuous and locally precompact relative to $D_P$, and $D_P \subset \dom{G_P}$;
        \item \label{item:4-prop} $\mathbb{U}_C$ and $\mathbb{U}_D$ are outer semicontinuous and locally precompact. 
    \end{enumerate}
    If $F_P\circ \mathbb{U}_C : \pi(C_P)\setto \T{}{\M}$ is convex valued, then $\hybrid_{\mathbb{U}}$ satisfies geometric hybrid basic conditions. 
\end{proposition}

\begin{proof}
    As $\mathcal{B}$ is a relaxed $C^1$-control bundle, $\M$ is a $C^1$-manifold and item~\ref{item:HBC1} of Definition~\ref{def:HBCs} holds. By Lemma~\ref{lem:mCBFFullDomains}, $C_{\mathbb{U}} = \pi(C_P)$ and $D_{\mathbb{U}} = \pi(D_P)$ are closed, so item~\ref{item:HBC2} holds. To verify that $F_{\mathbb{U}}$ satisfies item~\ref{HBC:F} of Definition~\ref{def:HBCs}, define $\widetilde{F}_{\mathbb{U}} : \M \setto \T{}{\M}$ as
    \begin{align*}
        \widetilde{F}_{\mathbb{U}}(x) \coloneqq \left\{
            \begin{array}{cc}
                F_{\mathbb{U}}(x) & \quad \text{if } x \in C_{\mathbb{U}}
                \\
                \varnothing & \quad \text{otherwise}
            \end{array}
        \right. \qquad \forall x\in \M. 
    \end{align*}
    Similarly, we also construct the map $\widetilde{F}_P : \U \setto \T{}{\M}$, defined for each $\chi\in \U$ as $\widetilde{F}_P(\chi) \coloneqq F_P(\chi)$ if $\chi \in C_P$, and $\widetilde{F}_P(\chi) = \varnothing$ otherwise. Then, by assumptions, $\widetilde{F}_P$ is outer semicontinuous and locally precompact at each $\chi\in C_P$.

    As $F_{\mathbb{U}}$ can be written as $F_{\mathbb{U}} = F_P\circ \mathbb{U}_C$, and as $\rge{\mathbb{U}_C} \subset C_P \subset \dom{F}$, it follows that
    \begin{align*}
        \widetilde{F}_{\mathbb{U}}(x) = (\widetilde{F}_P \circ {\mathbb{U}}_C)(x) = \left\{
            \begin{array}{cc}
                (F_P \circ \mathbb{U}_C)(x) & \quad \text{if } x\in C_{\mathbb{U}}
                \\
                \varnothing & \quad \text{otherwise}
            \end{array}
        \right. \qquad \forall x\in \M. 
    \end{align*}
    Consequently, using properties of $\widetilde{F}_{P}$ and ${\mathbb{U}}_C$ together with \cite[Thm. 17.23]{AliprantisBorder2006} and \cite[Prop. 6.3.2]{Beer1993}, it follows that $\widetilde{F}_P\circ \mathbb{U}_C$ is also outer semicontinuous and locally precompact at each $x\in \dom{\mathbb{U}_C}$. Therefore, $F_{\mathbb{U}}$ is outer semicontinuous and locally precompact. Finally, convex valuedness of $F_{\mathbb{U}}$ follows by assumption. Therefore, item~\ref{HBC:F} of Definition~\ref{def:HBCs} holds.  That item~\ref{HBC:G} of Definition~\ref{def:HBCs} holds follows similarly using properties of $G_{\mathbb{U}}$ and $D_{\mathbb{U}}$. 
\end{proof}

\begin{example}[Pendulum with impacts, revisited]
    \label{ex:pendulum-mCBF}
    Recall the open geometric hybrid system $\hybrid_P = (C_P, F_P, D_P, G_P, {h_P}, \mathcal{B})$ from Example~\ref{ex:pendulum-plant-dynamics}, and the mCBF candidate $B : \M \to \R{2}$ from Example~\ref{ex:pendulum-mCBF-candidate}. We construct the performance functions $\gamma = (\gamma_C, \gamma_D)$, $\gamma_C : \pi(C_P) \to \R{2}$ and $\gamma_D : \pi(D_P) \to \R{2}$, for which the function $B$ from Example~\ref{ex:pendulum-mCBF-candidate} is a mCBF for $\hybrid_P$ with respect to $(K, \gamma)$.  
    To this end, consider the saturation function $\sigma(s) \coloneqq \min\left\{1,{\max\{s,0\}}\right\}$ for each $s\in \R{}$. Define $\gamma_C:\pi(C_P)\to\R{2}$, where $\gamma_C = (\gamma_C^{\text{\tiny (-1)}}, \gamma_C^{\text{\tiny (1)}})$, as $\gamma_C^{\text{\tiny(-1)}}(x) \coloneqq -\sigma\bigl(E(x)-c_{-1}\bigr)$ and $\gamma_C^{\text{\tiny(1)}}(x) \coloneqq -\sigma\bigl(c_1-E(x)\bigr)$ for each $x\in \pi(C_P)$. Define also $\gamma_D:\pi(D_P)\to\R{2}$ as $\gamma_D^{\text{\tiny(i)}}(x)\coloneqq-B_i(x)$ for each $x\in \pi(D_P)$ and each $i\in\mathcal I = \{-1,1\}$. Let $\delta\in(0,(c_1-c_{-1})/2)$ and, for each $i\in\mathcal I$, define $M_i\coloneqq \{x\in\partial K:B_i(x)=0\}=\{x\in\M:E(x)=c_i\}$ and $U_i\coloneqq \{x\in\M: c_i - \delta < E(x) <  c_i+\delta\}$. 

    We verify that Assumption~\ref{ass:perfomanceFunctions} holds. For $i = -1$, pick any $x\in (U_{-1}\setminus K_{-1})\cap \pi(C_P)$. As $c_{-1} + \delta < c_1$, it follows that $E(x) <c_{-1}$, causing $\gamma_C^{\text{\tiny(-1)}}(x) = 0$. For $i = 1$, pick any $x\in (U_1\setminus K_1)\cap \pi(C_P)$. As $E(x) > c_1$, $\gamma_C^{\text{\tiny(1)}}(x) = 0$. Therefore, for each $i\in \mathcal{I}$ and each $x\in (U_i \setminus K_i)\cap \pi(C_P)$, $\gamma_C^{\text{\tiny(i)}}(x) = 0$. Similarly, for jumps, pick any $i\in \mathcal{I}$. For each $x\in K_{i} \cap \pi(D_P)$, $\gamma_D^{\text{\tiny(i)}}(x) = - B_{i}(x) \geq 0$ by construction of $K_i$. Therefore, Assumption~\ref{ass:perfomanceFunctions} holds. 

    To verify that $B$ is a mCBF for $\hybrid_P$ with respect to $(K, \gamma)$, note by construction that $C_P \subset \U = \dom{F_P}$ and $D_P \subset \U = \dom{G_P}$. Therefore, item~\ref{item:1-mCBF} of Definition~\ref{def:mCBF} holds.

    To verify item~\ref{item:2-mCBF}, pick any $x\in \pi(C_P)$ and let $q\coloneqq\pi_{\mathbb{S}^1}(x)$. As $x\in\T{q}{\mathbb{S}^1}$ is the velocity of the pendulum, choose $\chi_C\in\pi^{-1}(x)\cap C_P$ as
    \begin{align}
        \chi_C \coloneqq \left(x,bx-\grad{U}{q},0\right),
        \label{eq:pendulumFlowInput}
    \end{align}
    where $bx-\grad{U}{q}\in\T{q}{\mathbb{S}^1}$, and the gradient is computed with respect to the Riemannian metric $g_{\mathbb{S}^1}$. Substituting~\eqref{eq:pendulumFlowInput} in the definition of $F_P$ from Example~\ref{ex:pendulum-plant-dynamics} yields
    \[
        F_P(\chi_C) = \mathfrak{S}(x) + \vlift{x}{-\grad{U}{q}},
    \]
    where, recall that, $x = \pi(\chi_C)$ and $q = \pi_{\mathbb{S}^1}(x)$. Following the analysis in~\cite[Section~VIII]{JirwankarOchoaSanfelice_CDC2025} and using smoothness of $E$, we obtain that $\diffFunc{E}_{x}(F_P(\chi_C)) = 0$. Then, for each $i\in \mathcal{I}$,
    \begin{align*}
        \dot{B}_i(x) = \diffFunc{(B_i)}_{x}(F_P(\chi_C)) = i\; \diffFunc{E}_{x}(F_P(\chi_C)) = 0.
    \end{align*}
    As, for each $i\in \mathcal{I}$, $x\in (U_i\setminus K_i)\cap \pi(C_P)$ implies that $\gamma_C^{\text{\tiny(i)}}(x) = 0$, it follows that $\dot{B}_i(x) = 0 \leq -\gamma_C^{\text{\tiny(i)}}(x)$ for each $x\in (U_i\setminus K_i)\cap \pi(C_P)$ and each $i\in\mathcal{I}$. Alternatively, for each $i\in\mathcal{I}$ and each $x\in K_i \cap \pi(C_P)$, as $\gamma_C^{\text{\tiny(i)}}(x) \leq 0$ by construction, it follows that $\dot{B}_i(x) = 0 \leq -\gamma_C^{\text{\tiny(i)}}(x)$. Therefore, the inequality in item~\ref{item:2-mCBF} of Definition~\ref{def:mCBF} holds.

    To verify item~\ref{item:3-mCBF} of Definition~\ref{def:mCBF}, pick any $x \in \pi(D_P)$ and let $q\coloneqq\pi_{\mathbb{S}^1}(x)$. As the definition of $D_P$ does not constrain either control component, let
    \[
        \chi_D \coloneqq (x,0_q,1)\in \pi^{-1}(x)\cap D_P,
    \]
    where $0_q$ denotes the zero vector in $\T{q}{\mathbb{S}^1}$. Then, $G_P(\chi_D) = -\pi(\chi_D)=-x \in \T{}{\mathbb{S}^1}$. As the Riemannian inner product is a bilinear operator, $g_{\mathbb{S}^1}^{[z_1]}(v, v) = g_{\mathbb{S}^1}^{[z_1]}(-v, -v)$ for each $v\in \T{}{\mathbb{S}^1}$, where $[z_1] = \pi_{\mathbb{S}^1}(v)$. Therefore, using $v = \pi(\chi_D)=x$ in the Riemannian inner product, the total energy $E$ is preserved at jumps, causing $E(x) = E(G_P(\chi_D))$. Consequently, for each $i\in \mathcal{I}$,
    \begin{align}
        \sup_{g\in G_P(\chi_D)}B_i(g) = B_i(-x) = B_i(x) = -\gamma_D^{(i)}(x).
    \end{align}
    Therefore, $B$ is a mCBF for $\hybrid_P$ with respect to $(K, \gamma)$. Finally, to show that $G_P(D_P)\subset \pi(C_P) \cup \pi(D_P)$, pick any $\chi \in D_P$ and note that $G_P(\chi) = -\operatorname{pr}_3(\chi)\pi(\chi)$. Since $\operatorname{pr}_3(\chi)\in [0,1]$, scalar multiplication leaves the position of the pendulum unchanged. Moreover, $\pi_{\mathbb{S}^1}(\pi(\chi))\in\partial\cal W\subset\overline{\mathbb{S}^1\setminus\cal W}$. It follows from the definition of $C_P$ that $G_P(\chi)\in\pi(C_P)$. As this inclusion holds for each $\chi\in D_P$, we obtain that $G_P(D_P) \subset \pi(C_P)\subset \pi(C_P)\cup \pi(D_P)$. Then, using Lemma~\ref{lem:plantPostJumpCondition} with Theorem~\ref{thm:mCBFForwardPreInvariant}, $K$ is forward pre-invariant for $\hybrid_{\mathbb{U}}$.
\end{example}

\section{Pointwise Optimal and Safe Control Laws for $\hybrid_P$}
\label{sec:pointwiseOptimal}

Lemma~\ref{lem:singleSelectionInvariance} asserts that any control law that is a selection from $\Gamma$ renders $K$ forward pre-invariant for the resulting hybrid closed-loop system. However, the choice of a particular selection remains noncanonical. In this section, we make this choice canonical by picking only those selections that minimize a given cost function. In particular, given $Q : \U \to \R{}$, consider the set-valued map $\kappa : \dom{\kappa} \setto \U$ defined as
\begin{align}
    \label{eq:opti}
    \kappa(x) \coloneqq  \underset{\chi\in \Gamma(x)}{\arg\min} \; Q(\chi) \qquad \forall x\in \dom{\kappa}. 
\end{align}

From the construction of $\Gamma$ in~\eqref{eq:Gamma} and using the fact that $\mathbb{U}_C$ and $\mathbb{U}_D$ are nonempty valued on $\pi(C_P)$ and $\pi(D_P)$, respectively, it follows that the optimization problem in~\eqref{eq:opti} is feasible for each $x\in (\pi(C_P) \setminus \pi(C_P)) \cup (\pi(D_P)\setminus \pi(C_P))$. For points in $x\in \pi(C_P) \cap \pi(D_P)$, however, the intersection $\mathbb{U}_C(x) \cap \mathbb{U}_D(x)$ may be empty. Then, whether~\eqref{eq:opti} admits a minimum depends on properties of the cost $Q$ and the map $\Gamma$. To this end, we begin by stating the following assumption that, under additional conditions, asserts that the cost $Q$ is continuous and the aforementioned intersection is nonempty.

\begin{assumption}[Regularity of the pointwise optimization problems]
    \label{ass:measurableOptimalSelections}
    For the map $\Gamma : \pi(C_P)\cup\pi(D_P) \setto \U$ defined in~\eqref{eq:Gamma}, the set $\dom{\Gamma}$ is endowed with the trace $\sigma$-algebra
    \begin{align}
        \mathfrak{B}(\dom{\Gamma})
        &\coloneqq\{\dom{\Gamma}\cap A:A\in\mathfrak{B}(\M)\}, \label{eq:trace}
    \end{align}
    and the following conditions hold:
    \begin{enumerate}
        \item \label{item:optimalSelectionPolish} $\U$ is a Polish space;\footnote{A Polish space is a separable and completely metrizable topological space. A topological space is separable if it has a countable dense subset, and it is completely metrizable if there exists a complete metric on the space that is compatible with the topology.}
        \item \label{item:optimalSelectionCompact} $\Gamma(x)$ is a nonempty compact subset of $\U$ for each $x\in\pi(C_P) \cup \pi(D_P)$;
        \item \label{item:optimalSelectionMeasurable} $\Gamma$ is weakly measurable;\footnote{Let $(X, \mathfrak{A})$ be a measurable space and $Y$ a topological space. A map $F : X\rightrightarrows Y$ is weakly measurable if, for each open set $O \subset Y$, $F^{-1}(O) \subset \mathfrak{A}$~\cite[Definition~18.1]{AliprantisBorder2006}.} 
        \item \label{item:optimalSelectionContinuousCost} the function $Q:\U\to\R{}$ is continuous.
    \end{enumerate}
\end{assumption}

\begin{remark}
    Item~\ref{item:optimalSelectionPolish} of Assumption~\ref{ass:measurableOptimalSelections} is not restrictive. Indeed, that $\U$ is a Polish space holds for a wide class of spaces and, in particular, when $\U$ is a topological manifold. Therefore, it holds for Examples~\ref{ex:fiber-bundle},~\ref{ex:integrator},~\ref{ex:attitude-kinematics}, and~\ref{ex:pendulum-plant-dynamics}. 
\end{remark}

The above assumptions enable application of the measurable maximum theorem~\cite[Theorem~18.19]{AliprantisBorder2006} to show the existence of a measurable selection from $\kappa$.

\begin{theorem}[Measurable pointwise optimal  selections]
    \label{thm:measurableOptimalSelections}
    {Given a relaxed $C^1$-control bundle $\mathcal{B}=(\M,\U,\pi)$, an {open geometric hybrid system} $\hybrid_P=(C_P,F_P,D_P,G_P,{h_P},\mathcal{B})$, a nonempty closed set $K \subset \M$, performance functions $\gamma=(\gamma_C,\gamma_D)$ satisfying Assumption~\ref{ass:perfomanceFunctions}, a mCBF $B:\dom{B}\to\R{N}$ for $\hybrid_P$ with respect to $(K,\gamma)$, and the set-valued map $\Gamma : \pi(C_P)\cup \pi(D_P) \setto \U$ as in~\eqref{eq:Gamma}, suppose that Assumption~\ref{ass:measurableOptimalSelections} holds. Then, $\dom{\Gamma} =\pi(C_P) \cup \pi(D_P)$, and}
    \begin{enumerate}[label=(\roman*)]
        \item \label{item:thm-kappa-weaklyMeasurable-compact} {the map $\kappa$ in~\eqref{eq:opti}} is such that $\dom{\kappa} = \dom{\Gamma}$, $\kappa(x)$ is nonempty and compact for each $x\in \dom{\kappa}$, and $\kappa$ is weakly measurable;
        \item there exist Borel-measurable functions $\widehat{\kappa}:\dom{\Gamma}\to\U$ such that
    \begin{align}
        \label{eq:measurableOptimalSelection}
        \widehat{\kappa}(x)&\in\kappa(x) \qquad \forall x\in\dom{\Gamma}.
    \end{align}
    \end{enumerate}
\end{theorem}

\begin{proof}
    By item~\ref{item:optimalSelectionCompact} of Assumption~\ref{ass:measurableOptimalSelections}, $\dom{\Gamma} = \pi(C_P) \cup \pi(D_P)$. Let $x\in\dom{\Gamma}$ be arbitrary. Item~\ref{item:optimalSelectionCompact} of Assumption~\ref{ass:measurableOptimalSelections} implies that $\Gamma(x)$ is nonempty and compact. Item~\ref{item:optimalSelectionContinuousCost} of Assumption~\ref{ass:measurableOptimalSelections} implies that $Q|_{\Gamma(x)}$ is continuous. {Then, the extreme value theorem~\cite[Theorem~2.43]{AliprantisBorder2006} yields $\kappa(x)\neq\varnothing$. As $\kappa(x) \subset \Gamma(x)$ by construction, and $\kappa(x)$ collects the minimizers of $Q$ on $\Gamma(x)$, it follows that 
    \begin{align*}
        \kappa(x)
        =\Gamma(x)\cap
        Q^{-1}(q)
    \end{align*}
    where $q \coloneqq \min_{v\in \Gamma(x)}Q(v) \in \R{}$ denotes the minimum value of $Q$ on $\Gamma(x)$. As $q\in\R{}$ is closed and $Q$ is continuous by item~\ref{item:optimalSelectionContinuousCost} of Assumption~\ref{ass:measurableOptimalSelections}, the preimage $Q^{-1}(q)$ is closed in $\U$.} Next, as $\Gamma(x)$ is compact and $\U$ is Hausdorff as it is a Polish space, it follows from~\cite[Proposition~A.45(f)]{Lee} that $\Gamma(x)$ is closed. Consequently, $\kappa(x)$ is a closed subset of the compact set $\Gamma(x)$. As $x\in\dom{\Gamma}$ was arbitrary, $\kappa$ has nonempty, closed, and compact values on $\dom{\Gamma}$. Consequently, $\dom{\kappa} = \dom{\Gamma}$ and $\kappa(x)$ is nonempty and compact for each $x\in\dom{\Gamma}$. 

    Define $\widetilde{q}:\graph{\Gamma}\to\R{}$ as $\widetilde{q}(x,u)\coloneqq-Q(u)$ for each $(x,u)\in\graph{\Gamma}$. To apply the measurable maximum theorem~\cite[Theorem~18.19]{AliprantisBorder2006}, extend $\widetilde{q}$ to the product space $\dom{\Gamma}\times \U$ by defining $q:\dom{\Gamma}\times\U\to\R{}$ as $q(x,u)\coloneqq-Q(u)$ for each $x\in \dom{\Gamma}$ and each $u\in \U$. Note that, for arbitrary choices of $x\in \dom{\Gamma}$ and $u\in \U$, it may so happen that $u\notin \pi^{-1}(x)$. For each $u\in\U$, the function $x\mapsto q(x,u)$ is constant and therefore $(\mathfrak{B}(\dom{\Gamma}), \mathfrak{L}(\R{}))$-measurable. Then, using $\mathfrak{B}(\R{})\subset \mathfrak{L}(\R{})$, $x\mapsto q(x,u)$ is $(\mathfrak{B}(\dom{\Gamma}),\mathfrak{B}(\R{}))$-measurable for each $u\in \U$. For each $x\in\dom{\Gamma}$, item~\ref{item:optimalSelectionContinuousCost} of Assumption~\ref{ass:measurableOptimalSelections} implies that the function $u\mapsto q(x,u)$ is continuous. Therefore, $q$ is a Carath\'eodory function~\cite[Definition~4.50]{AliprantisBorder2006}.\footnote{Let $(S, \Sigma)$ be a measurable space and $X$ and $Y$ be topological spaces. A function $f : S\times X \to Y$ is a Carath\'eodory function if (i) for each $x\in X$, the function $s\mapsto f(s,x)$ is $(\Sigma, \mathfrak{B}(Y))$-measurable, and (ii) for each $s\in S$, the function $x\mapsto f(s,x)$ is continuous.} Item~\ref{item:optimalSelectionMeasurable} of Assumption~\ref{ass:measurableOptimalSelections} states that $\Gamma$ is measurable, and item~\ref{item:optimalSelectionCompact} states that its values are nonempty and compact. Moreover, item~\ref{item:optimalSelectionPolish} implies that $\U$ is a separable metrizable space. The measurable maximum theorem applied to $q$ and $\Gamma$ yields that
    \begin{align}
        x\mapsto\underset{u\in\Gamma(x)}{\arg\max}\;q(x,u)
        =\kappa(x)
    \end{align}
    is a nonempty, compact-valued measurable map. Then, by~\cite[Lemma~18.2]{AliprantisBorder2006}, measurability of $\kappa$ also implies that it is weakly measurable. Additionally, item~(iii) of the measurable maximum theorem~\cite[Theorem~18.19]{AliprantisBorder2006} concludes that $\kappa$ admits a measurable selection $\widehat{\kappa}$. 

\end{proof}

Items~\ref{item:thm-kappa-weaklyMeasurable-compact} of Theorem~\ref{thm:measurableOptimalSelections} yields that $\kappa$ is a weakly measurable map with nonempty, closed values. These properties of $\kappa$ are crucial for the application of the measurable maximum theorem that certifies the existence of measurable selections from $\kappa$. It is interesting to note that, if $\kappa$ is outer semicontinuous and $\U$ is $\sigma$-compact,\footnote{A topological space is $\sigma$-compact if it is the union of countably many compact sets.} then it is weakly measurable and closed-valued~\cite[Theorem~18.20]{AliprantisBorder2006}, enabling the application of the selection theorem. Therefore, we present the following conditions under which $\kappa$ is outer semicontinuous. 

\begin{proposition}
    \label{prop:osc}
    {Given a relaxed $C^1$-control bundle $\mathcal{B}=(\M,\U,\pi)$, an {open geometric hybrid system} $\hybrid_P=(C_P,F_P,D_P,G_P,{h_P},\mathcal{B})$, a nonempty closed set $K \subset \M$, performance functions $\gamma=(\gamma_C,\gamma_D)$ satisfying Assumption~\ref{ass:perfomanceFunctions}, a mCBF $B:\dom{B}\to\R{N}$ for $\hybrid_P$ with respect to $(K,\gamma)$, and the set-valued map $\Gamma : \pi(C_P)\cup \pi(D_P) \setto \U$ as in~\eqref{eq:Gamma}, suppose that $\U$ is Hausdorff and $\dom{\Gamma} =\pi(C_P) \cup \pi(D_P)$.} Suppose also that 
    \begin{enumerate}[label=(\roman*)]
        \item $\Gamma$ is continuous and locally precompact, and \item $Q$ is continuous.
    \end{enumerate}
    Then, $\dom{\kappa} = \dom{\Gamma}$ and $\kappa$ is outer semicontinuous and locally precompact.
\end{proposition}

\begin{proof}
    By the domain equality and using outer semicontinuity, $\Gamma : \pi(C_P) \cup \pi(D_P) \setto \U$ has nonempty closed values everywhere on $\pi(C_P) \cup \pi(D_P)$. Additionally, using local precompactness, the closed values of $\Gamma$ are subsets of some compact set. Then, by~\cite[Proposition~A.45(e)]{Lee}, $\Gamma$ is also compact valued. Finally, outer semicontinuity and local precompactness of $\Gamma$ also cause $\Gamma$ to be upper continuous. Together with inner semicontinuity, these maps are continuous in the sense of~\cite[Definition~17.2]{AliprantisBorder2006}. By Berge's maximum principle~\cite[Theorem~17.31]{AliprantisBorder2006}, $\kappa$ is upper hemicontinuous;  see~\cite[Definition~17.2]{AliprantisBorder2006}. As $\U$ is Hausdorff and $\Gamma$ is compact valued, an application of~\cite[Theorem~17.10]{AliprantisBorder2006} causes $\kappa$ to be outer semicontinuous. Local precompactness of $\kappa$ follows from local precompactness of $\Gamma$ by noting that $\kappa(x)\subset \Gamma(x)$ for each $x\in \dom{\Gamma}$. Finally, that $\kappa$ is nonempty valued; namely, $\dom{\kappa} = \dom{\Gamma} = \pi(C_P) \cup \pi(D_P)$, follows from the extreme value theorem as $Q$ is continuous and $\Gamma$ is compact valued, causing the minima in~\eqref{eq:opti} to exist for each $x\in \dom{\Gamma}$. 
\end{proof}

Note that Theorem~\ref{thm:measurableOptimalSelections} requires $\U$ to be a Polish space, but Proposition~\ref{prop:osc} requires $\U$ to be only Hausdorff. Indeed, as each Polish space is metrizable by definition, it is also Hausdorff~\cite[\S21, p.~129]{Munkres2000}. Additionally, under $\sigma$-compactness of $\U$, it immediately follows from Proposition~\ref{prop:osc} that there exist measurable selections from $\kappa$. Note that $\sigma$-compactness of $\U$ is not restrictive, and it holds in particular when $\U$ is a topological manifold. 

\begin{corollary}
    \label{cor:selection}
    {Given a relaxed $C^1$-control bundle $\mathcal{B}=(\M,\U,\pi)$, an {open geometric hybrid system} $\hybrid_P=(C_P,F_P,D_P,G_P,{h_P},\mathcal{B})$, a nonempty closed set $K \subset \M$, performance functions $\gamma=(\gamma_C,\gamma_D)$ satisfying Assumption~\ref{ass:perfomanceFunctions}, a mCBF $B:\dom{B}\to\R{N}$ for $\hybrid_P$ with respect to $(K,\gamma)$, and the set-valued map $\Gamma : \pi(C_P)\cup \pi(D_P) \setto \U$ as in~\eqref{eq:Gamma}, suppose that $\U$ is a $\sigma$-compact Polish space, $\dom{\Gamma} =\pi(C_P) \cup \pi(D_P)$, $\Gamma$ is a continuous and locally precompact set-valued map, $Q$ is continuous, and $\dom{\Gamma}$ is endowed with the trace $\sigma$-algebras as in~\eqref{eq:trace}.} Then, there exists a measurable selection from $\kappa$.
\end{corollary}

\begin{proof}
    By Proposition~\ref{prop:osc}, $\kappa : \pi(C_P)\cup\pi(D_P) \setto \U$ satisfies $\dom{\kappa} = \pi(C_P)\cup\pi(D_P)$, $\kappa$ is outer semicontinuous, and has nonempty closed values. By~\cite[Theorem~18.20]{AliprantisBorder2006}, they are also weakly measurable. Then, the Kuratowski--Ryll-Nardzewski selection theorem~\cite[Theorem~18.13]{AliprantisBorder2006} yields a measurable selection from $\kappa$. 
\end{proof}

\begin{remark}[Relation to the Euclidean case]
    Note that the {state-input point} in~\eqref{eq:hybridplant} is denoted by $\chi\in \U$. In the Euclidean setting~\cite{montenegro2025pointwise}, such a $\chi$ is factorized into flow and jump input components. Then, assigning the flow input from $\mathbb{U}_C$ and the jump input from $\mathbb{U}_D$, and requiring Assumption~\ref{ass:measurableOptimalSelections} to hold for $\mathbb{U}_C$ and $\mathbb{U}_D$, one can obtain the above certificates without requiring $\Gamma$ in~\eqref{eq:Gamma}. Therefore, comparing the above results with~\cite{montenegro2025pointwise}, we observe that Assumption~\ref{ass:measurableOptimalSelections} and Theorem~\ref{thm:measurableOptimalSelections} generalize the Euclidean results~\cite{montenegro2025pointwise} to obtain a measurable selections instead of a continuous one as in~\cite{montenegro2025pointwise}. Additionally, Proposition~\ref{prop:osc} extends~\cite[Lemma~4.2]{montenegro2025pointwise} to the geometric setting. The stronger assumption in~\cite[Lemma~4.2]{montenegro2025pointwise} on strict convexity of $Q$ the image of $\mathbb{U}_C$ and $\mathbb{U}_D$, respectively yield single valued and continuous flow and jump control laws. While the notion of convexity is not suitable for a geometric analysis as it requires the underlying space to be linear, enforcing that~\eqref{eq:opti} has unique optimizers for each $x$ in can lead to a result similar to~\cite[Lemma~4.2]{montenegro2025pointwise}. However, if $\U$ has a vector bundle structure, each fiber of $\pi$ is a vector space, and convexity of $Q$ may be enforced on each fiber to yield a continuous selection. In the Euclidean setting, the such continuity of $\kappa$ enables convenient certificates for robustness to structural perturbations~\cite[Ch.~6]{goebel_hybrid_2012}. Such a robustness framework for geometric hybrid systems is under development and, therefore, the existence of continuous selections is part of the future work. 
\end{remark}

\begin{remark}[Relation to CBF-QP problems]
    In the Euclidean setting for continuous-time differential equations with control-affine plant dynamics, one has $\dot{z} = f(z) + g(z)u$, where $z\in \R{n}$ denotes the state and the input $u\in \R{m}$ is assumed to be unconstrained. Then, with a cost that is quadratic in $u$, the optimization problem~\eqref{eq:opti} can be solved easily -- even in real-time for onboard implementation --  and yields a continuous controller that, consequently, provides desirable robustness guarantees. Similar results in the geometric setting are found in~\cite{ames2025bundlesbacksteppinggeometriccontrol}. They can be extended to open geometric hybrid systems on control bundles when $\U$ has a vector bundle structure (see Example~\ref{ex:vectorBundle}). This is possible because the optimization problem~\eqref{eq:opti} is solved on a fiber of $\pi$ which, under the vector bundle structure, becomes a vector space. 
\end{remark}

\section{Application to Geometric Hybrid Safety Filters}
\label{sec:safetyFilter}

Recall from the previous section that the control law $\widehat{\kappa}$, which is a selection from $\kappa$, was obtained by minimizing the cost function $Q : \U \to \R{}$. This cost function can be designed to encode information about nominal static state-feedback control law $\mu :\M \to \U$ as in Definition~\ref{def:staticFeedback} so that $\widehat{\kappa}$ minimally modifies the nominal controller $\mu$ to yield a \emph{filtered} controller that renders the set $K$ forward invariant for the resulting closed-loop system. The above setup, called a \emph{safety filter}, has rapidly gained popularity as a tool to leverage uncertified controllers with better performance (for example, those obtained via reinforcement learning~\cite{Schoellig2025SafetyRL,choi2020reinforcement}, imitation learning~\cite{Ames2022imitation}, or arbitrary neural-network-based controllers~\cite{Fleming2025ShieldNN}) while also ensuring safety via the optimization problems~\eqref{eq:opti} through mCBF decrease conditions. 

Next, we formalize the above setup. Given a static state-feedback control law $\mu : \M \to \U$ that is not necessarily safeguarding, consider the cost functions $Q_{\mu} : \U \to \R{}_{\geq 0}$ that satisfies the following assumption.

\begin{assumption}[Safety-filter cost]
    \label{ass:safetyFilterCost}
    Consider a relaxed $C^1$-control bundle $(\M, \U, \pi)$. Then, a static state-feedback control law $\mu : \M \to \U$ and the cost function $Q_{\mu} : \U \to \R{}_{\geq 0}$ satisfy, for each ${\chi} \in \U$ and letting $x\coloneqq \pi({\chi})$, that
    
    \begin{align}
        \label{eq:filterCost}
        Q_{\mu}(\chi) = 0 \quad \text{if and only if} \quad \chi = \mu(x). 
    \end{align}
\end{assumption}

The subscript $\mu$ in $Q_{\mu}$ highlights the dependence of the cost on the nominal controller $\mu$. The condition in~\eqref{eq:filterCost} ensures that the cost is zero if and only if its argument equals the corresponding nominal control input. Therefore, the above assumption enforces positive definiteness of $Q_{\mu}$ with respect to $\rge{\mu}$. 

\begin{lemma}[Positive definiteness of the safety-filter cost]
    Given a relaxed $C^1$-control bundle $(\M, \U, \pi)$, a static state-feedback control law $\mu : \M \to \U$, and a cost function $Q_{\mu} : \U \to \R{}_{\geq 0}$, suppose that Assumption~\ref{ass:safetyFilterCost} holds. Then, $Q_{\mu}$ is positive definite with respect to $\rge{\mu}$. 
\end{lemma}

\begin{proof}
    Pick any ${\chi} \in \U$. Suppose that $\chi \in \rge{\mu}$, and let $x \coloneqq \pi({\chi})$. By definition, $\chi = \mu(x)$. Then, by Assumption~\ref{ass:safetyFilterCost}, $Q_{\mu}(\chi) = 0$. Alternatively, suppose that $\chi \notin \rge{\mu}$. Then, note that $\chi\neq \mu(x)$ and, therefore, $\chi\neq \mu(\pi({\chi}))$. Then, Assumption~\ref{ass:safetyFilterCost} causes $Q_{\mu}(\chi) > 0$.
\end{proof}

Next, following the approach in Section~\ref{sec:pointwiseOptimal}, we define the set-valued map $\kappa_{\mathrm{sf}} : \dom{\Gamma} \rightrightarrows \U$ by
\begin{align}
    \kappa_{\mathrm{sf}}(x)
        &\coloneqq \underset{\chi\in\Gamma(x)}{\arg\min}\; Q_{\mu}(\chi)
        && \forall x\in\dom{\Gamma}.
        \label{eq:filter}
\end{align}

It is crucial to note that $Q_{\mu}$ captures the ``difference'' between a control input and the nominal controller $\mu$ in the same fiber. This notion of ``difference'' does not require a Riemannian metric or a vector space structure on $\U$. However, when such structure is available, it can be leveraged to define $Q_{\mu}$ suitably; see forthcoming Example~\ref{ex:pendulum-safety-filter}.

Under Assumption~\ref{ass:measurableOptimalSelections}, with $Q$ therein replaced by $Q_{\mu}$, the above map is nonempty-valued everywhere on its domains. By positive definiteness, if $\mu(x)\in\Gamma(x)$, then $\mu(x)\in\kappa_{\mathrm{sf}}(x)$. Therefore, wherever the corresponding nominal controller satisfies the mCBF constraints, the safety filter does not modify that nominal control input.

With the safety-filtered set-valued control law defined in~\eqref{eq:filter}, and under Assumption~\ref{ass:measurableOptimalSelections}, Theorem~\ref{thm:measurableOptimalSelections} yields a measurable selection $\widehat{\kappa}_{\mathrm{sf}}$ from $\kappa_{\mathrm{sf}}$ that, when applied to the hybrid plant $\hybrid_P$, yields the hybrid closed-loop system $\hybrid_{\mathrm{sf}} = (C_{\mathrm{sf}}, F_{\mathrm{sf}}, D_{\mathrm{sf}}, G_{\mathrm{sf}}, \M)$, where
\begin{align}
    \hybrid_{\mathrm{sf}} : \left\{
        \begin{array}{rcll}
            \dot{x} & \in F_{\mathrm{sf}}(x) \coloneqq F_P(\widehat{\kappa}_{\mathrm{sf}}(x))
                & \qquad & x\in {C_{\widehat{\kappa}_{\mathrm{sf}}},}
            \\[10pt]
            x^+ & \in G_{\mathrm{sf}}(x) \coloneqq G_P(\widehat{\kappa}_{\mathrm{sf}}(x))
                & \qquad & x\in {D_{\widehat{\kappa}_{\mathrm{sf}}}.}
        \end{array}
    \right.
    \label{eq:safetyFilteredClosedLoop}
\end{align}
where $C_{\widehat{\kappa}_{\mathrm{sf}}} \coloneqq \{x\in\pi(C_P)\cup \pi(D_P) : \widehat{\kappa}_{\mathrm{sf}}(x) \in C_P\}$ and $D_{\widehat{\kappa}_{\mathrm{sf}}} \coloneqq \{x\in\pi(C_P)\cup \pi(D_P) : \widehat{\kappa}_{\mathrm{sf}}(x) \in D_P\}$. The following result establishes forward pre-invariance of $K$ for $\hybrid_{\mathrm{sf}}$. 
\begin{theorem}[Forward pre-invariance under a safety-filtered selection]
    \label{thm:safetyFilteredForwardPreInvariance}
    Given a relaxed $C^1$-control bundle $\mathcal{B} = (\M, \U, \pi)$, an {open geometric hybrid system} $\hybrid_P = (C_P, F_P, D_P, G_P, {h_P}, \mathcal{B})$, a nonempty closed set $K\subset\M$, and performance functions $\gamma = (\gamma_C, \gamma_D)$ satisfying Assumption~\ref{ass:perfomanceFunctions}, let $B : \dom{B}\to\R{N}$ be a mCBF for $\hybrid_P$ with respect to $(K,\gamma)$. Let $\mu:\M\to\U$ be a static state-feedback controller, and let $Q_{\mu}:\U\to\R{}_{\geq0}$ satisfy Assumption~\ref{ass:safetyFilterCost}. Suppose also that Assumption~\ref{ass:measurableOptimalSelections} holds with $Q$ therein replaced by $Q_{\mu}$. If
    \begin{align}
        \label{eq:jumpInclusion}
        G_P(D_P) \subset \pi(C_P) \cup \pi(D_P),
    \end{align}
    then there exists a Borel-measurable selection $\widehat{\kappa}_{\mathrm{sf}}$ from $\kappa_{\mathrm{sf}}$. Moreover, for each such selection, $K$ is forward pre-invariant for the resulting hybrid closed-loop system $\hybrid_{\mathrm{sf}}$ in~\eqref{eq:safetyFilteredClosedLoop}. 
\end{theorem}

\begin{proof}
    Theorem~\ref{thm:measurableOptimalSelections}, applied to the map
    $\kappa_{\mathrm{sf}}$ in~\eqref{eq:filter}, yields a Borel-measurable
    selection $\widehat{\kappa}_{\mathrm{sf}}:\dom{\Gamma}\to\U$
    satisfying $\widehat{\kappa}_{\mathrm{sf}}(x) \subset \Gamma(x)$ for each $x\in\dom{\Gamma}$. 
    Item~\ref{item:optimalSelectionCompact} of
    Assumption~\ref{ass:measurableOptimalSelections} additionally yields $\dom{\Gamma}=\pi(C_P)\cup\pi(D_P)$. Consequently, we have that
    $\widehat{\kappa}_{\mathrm{sf}}$ is a selection from $\Gamma$
    on $\pi(C_P)\cup\pi(D_P)$.

    Next, Lemma~\ref{lem:plantPostJumpCondition} and~\eqref{eq:jumpInclusion}
    yield that $G_{\mathbb{U}}(K\cap D_{\mathbb{U}}) \subset C_{\mathbb{U}}\cup D_{\mathbb{U}}$. 
    Therefore, the assumptions of
    Theorem~\ref{thm:mCBFForwardPreInvariant} hold. Then, applying Lemma~\ref{lem:singleSelectionInvariance} to $\widehat{\kappa}_{\mathrm{sf}}$ yields $C_{\widehat{\kappa}_{\mathrm{sf}}}=\pi(C_P)$ and $D_{\widehat{\kappa}_{\mathrm{sf}}}=\pi(D_P)$, and it certifies that $K$ is forward pre-invariant for the closed-loop
    system $\hybrid_{\mathrm{sf}}$. 
\end{proof}

\begin{example}[Pendulum with impacts, revisited]
    \label{ex:pendulum-safety-filter}
    Recall the {open geometric hybrid system} $\hybrid_P$ in Example~\ref{ex:pendulum-plant-dynamics}. Consider the nominal controller $\mu : \M \to \U$ defined as 
    \begin{align}
        \label{eq:nominalController}
        \mu(z) \coloneqq \left(z,-\grad{U}{\pi_{\mathbb{S}^1}(z)},0.8\right)
        \qquad \forall z\in \M.
    \end{align}
    Then, note that $\operatorname{pr}_2(\mu(z)) \in \T{\pi_{\mathbb{S}^1}(z)}{\mathbb{S}^1}$ {denotes the force due to gravity obtained from the gravitational potential energy $U$}, and $\operatorname{pr}_3(\mu(z)) \in [0,1]$ denotes the controlled coefficient of restitution, which is set to constant $0.8$ in~\eqref{eq:nominalController}. {This choice of $\mu$ models the motion of the pendulum under gravity, with inelastic impacts at the wall.} The resulting closed-loop system is studied in~\cite[Section~VIII]{JirwankarOchoaSanfelice_CDC2025}, and solutions to it asymptotically approach $\A \coloneqq \{z\in \pi(C_P)\cup \pi(D_P) : E(z) = c_0\}$, where recall from Example~\ref{ex:pendulum-mCBF-candidate} that $c_0 = E(\partial \cal W \times \{0\})$ and $c_0 < c_{-1}$. Therefore, the above nominal control law $\mu$ is not safeguarding as solutions to its closed-loop system exit $K$ to approach the set $\A$. This fact is observed in Figure~\ref{fig:pendulum}, where the solutions using the nominal controller exit the safe set, colored in green. 

    Then, to construct the safety filter, define the cost $Q_{\mu} : \U \to \R{}_{\geq 0}$ as
    \begin{align}
        Q_{\mu}(\chi)
        &\coloneqq \left|\operatorname{pr}_2(\chi)-\operatorname{pr}_2(\mu(\pi(\chi)))\right|_{\mathbb{S}^1}^2
        +\left|\operatorname{pr}_3(\chi)-\operatorname{pr}_3(\mu(\pi(\chi)))\right|^2
        \qquad \forall \chi\in \U.
        \label{eq:cost-safetyFilter}
    \end{align}
    where $|v|^2_{\mathbb{S}^1} \coloneqq g_{\mathbb{S}^1}(v, v)$ for each $v\in\T{}{\mathbb{S}^1}$ denotes the vector norm in the tangent space using the Riemannian metric, and $|\cdot|^2$ denotes the Euclidean norm on $\R{}$. Therefore, $Q_\mu(\chi)$ captures the ``difference'' between $\chi$ and the nominal input $\mu(\pi(\chi))$ in the same fiber, where its first term measures the difference between force components of the nominal controller and the safety-filtered controller, while its second term measures the difference between their coefficients of restitution. It follows that, for each $\chi\in \U$, $Q_\mu(\chi)=0$ if and only if $\chi=\mu(\pi(\chi))$, and Assumption~\ref{ass:safetyFilterCost} holds.
    
    \begin{figure}
        \centering
        {\includegraphics[width=\linewidth, trim = 0cm 6.7cm 0cm 0cm, clip]{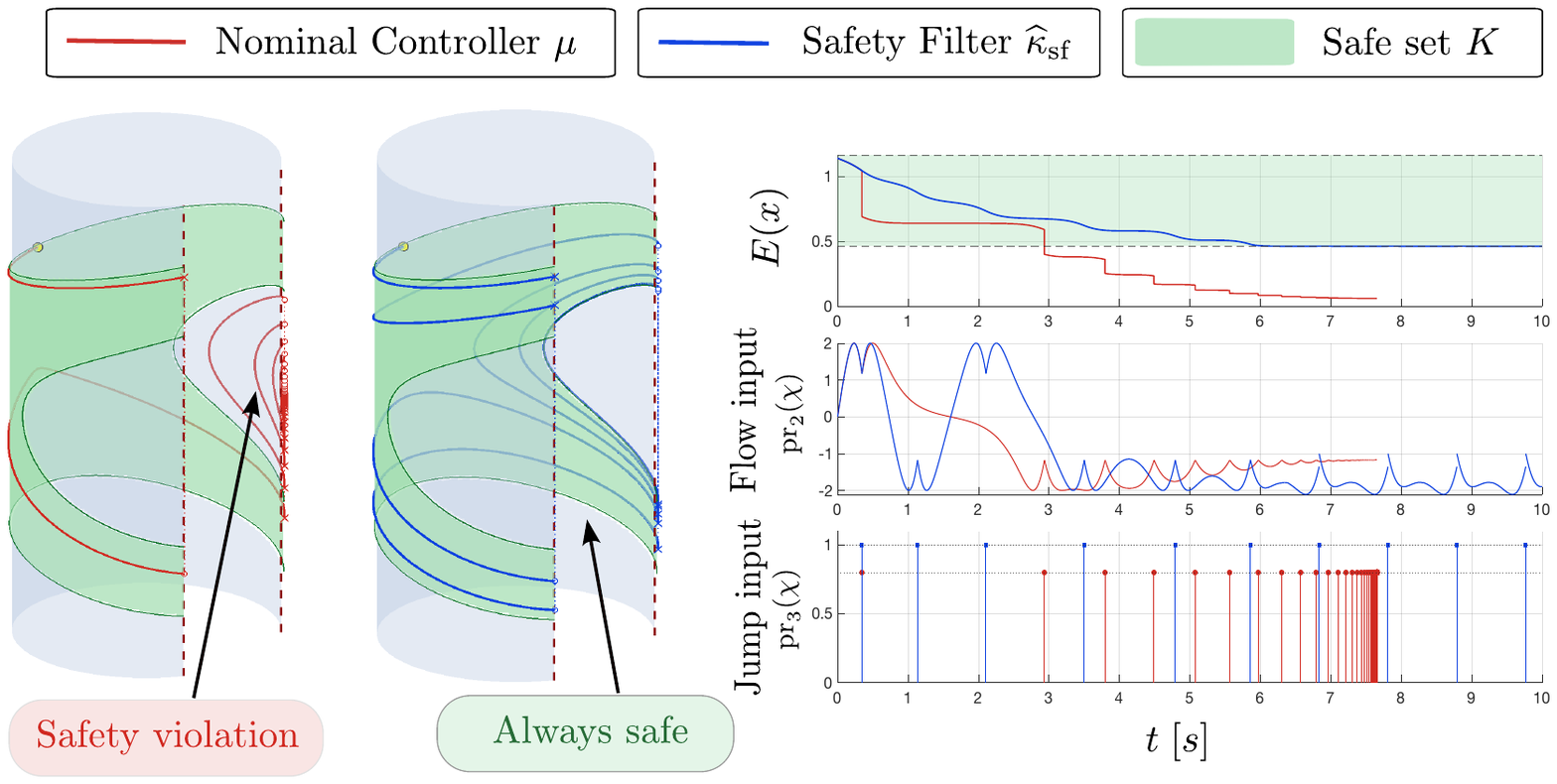}}
        \caption{
            \emph{Left:} $\M = \T{}{\mathbb{S}^1}$ denoted by the cylinder. The safe set $K$ is depicted in green. Solutions to $\hybrid_P$ using the nominal controller exit $K$. 
            \emph{Center:} Trajectories in blue, that denote solutions to $\hybrid_P$ using the safety-filtered control law $\widehat{\kappa}_{\mathrm{sf}}$, remain in $K$. 
            \emph{Right, top:} Total energy along solutions vs time. The nominal controller (red) causes the solution to exit the safe set $K$. The safety-filtered controller keeps the solutions inside $K$.
            \emph{Right, middle:} Flow input component vs time, where the former is projected onto $\R{}$ by identifying it with the tangent space to $\mathbb{S}^1$.
            \emph{Right, bottom:} Jump input component vs time. The nominal solution exhibits Zeno behavior as observed by increasing number of jumps in finite time, which is absent in the safety-filtered solution.
        }
        \label{fig:pendulum}
    \end{figure}
    
    As the projection operators and $\mu$ are continuous, and the Riemannian inner product is smooth, the cost function $Q_{\mu}$ is continuous, causing item~\ref{item:optimalSelectionContinuousCost} of Assumption~\ref{ass:measurableOptimalSelections} to hold. Additionally, as $\T{}{\mathbb{S}^1}\oplus\T{}{\mathbb{S}^1}$ is a smooth manifold~\cite[Example~10.7]{Lee}, the product $\U=(\T{}{\mathbb{S}^1}\oplus\T{}{\mathbb{S}^1})\times[0,1]$ is a Polish space. Therefore, item~\ref{item:optimalSelectionPolish} of Assumption~\ref{ass:measurableOptimalSelections} holds. Item~\ref{item:optimalSelectionContinuousCost} of Assumption~\ref{ass:measurableOptimalSelections} also holds as $Q_{\mu}$ is continuous. Under items~\ref{item:optimalSelectionCompact} and~\ref{item:optimalSelectionMeasurable} of Assumption~\ref{ass:measurableOptimalSelections}, and using $G_P(D_P)\subset \pi(C_P) \cup \pi(D_P)$ from Example~\ref{ex:pendulum-mCBF}, Theorem~\ref{thm:safetyFilteredForwardPreInvariance} yields the existence of a measurable selection $\widehat{\kappa}_{\mathrm{sf}}$ from the safety-filter map in~\eqref{eq:filter} that renders the set $K$ forward pre-invariant for the resulting closed-loop system. Figure~\ref{fig:pendulum} illustrates this result, with the various parameters set as $a \coloneqq 2$, $b \coloneqq 0.2$,  $c_0 = 0.06$, $c_{-1} \coloneqq 0.46$, and $c_{1} \coloneqq 1.16$.
\end{example}

\section{Conclusions and Future Work}

We have proposed a general framework to model geometric hybrid systems on manifolds with inputs. This framework extends control bundles to define \emph{relaxed} $C^1$-control bundles, which allow inputs to belong to an arbitrary topological space. This generality is crucial to model hybrid dynamics, where the input for flows and jumps may have different ``dimensions'' that can be conveniently captured by a topological space. Using such geometric hybrid systems with inputs, we have obtained multiple control barrier function (mCBF)-based, pointwise-optimal control laws, and conditions for existence of their Borel-measurable selections. To analyze the closed-loop system, we have also presented multiple barrier function-based sufficient conditions for forward pre-invariance of a closed set. Finally, we have applied the mCBF-based safe control design method to obtain safety filters for hybrid dynamical systems on manifolds. Together, these results present a systematic toolset for safe control synthesis for hybrid systems on manifolds.

Future directions involve obtaining sufficient conditions for existence of continuous selections of control laws, instead of merely Borel-measurable selections, and studying their effect on robustness to structural perturbations. The proposed framework with inputs also provides opportunities to investigate interconnections of hybrid systems on manifolds, which has applications in designing hybrid control laws instead of static state-feedback laws that are considered in this paper.

\backmatter

\bmhead{Supplementary information}

{We have included~\cite{jirwankar2026barrier} in the supplementary material. It has been accepted for publication at 2026 IEEE 65th Conference on Decision and Control (CDC).
}

\bmhead{Acknowledgements}

Research partially supported by NSF Grants no.\ CNS-2039054 and CNS-2111688, by AFOSR Grants nos.\ FA9550-19-1-0169, FA9550-20-1-0238, FA9550-23-1-0145, and FA9550-23-1-0313, by AFRL Grant nos.\ FA8651-22-1-0017 and FA8651-23-1-0004, by ARO Grant no.\ W911NF-20-1-0253, and by DoD Grant no.\ W911NF-23-1-0158.

\section*{Declarations}

\section*{Funding}
Research partially supported by NSF Grants no.\ CNS-2039054 and CNS-2111688, by AFOSR Grants nos.\ FA9550-19-1-0169, FA9550-20-1-0238, FA9550-23-1-0145, and FA9550-23-1-0313, by AFRL Grant nos.\ FA8651-22-1-0017 and FA8651-23-1-0004, by ARO Grant no.\ W911NF-20-1-0253, and by DoD Grant no.\ W911NF-23-1-0158.

\section*{Conflict of interest/Competing interests}
 The authors declare no conflict of interest.

\section*{Ethics approval and consent to participate}
Not applicable.

\section*{Consent for publication}
Not applicable.

\section*{Data availability }
Not applicable.

\section*{Materials availability}
Not applicable.

\section*{Code availability}
The simulation files are available at \url{https://github.com/pjirwank/safety-filter-pendulum-with-impacts}. 

\section*{Author contribution}
{P.P.J. contributed to the writing of the manuscript. P.P.J., D.E.O., C.A.M.G., and R.G.S. contributed to the mathematical analysis. All authors reviewed and approved the final version of the manuscript.}

\begin{appendices}

\section{Modeling Geometric Continuous-time and Discrete-time Systems}
\label{app:modeling-CT-DT}

Given a relaxed $C^1$-control bundle, the modeling of an open geometric hybrid system in~\eqref{eq:hybridplant} reduces to a differential inclusion (resp., difference inclusion) when the jump set $D_P$ (resp., the flow set $C_P$) is empty. In this section, we study such modeling by representing the state-input information as a tuple $(z,\chi)$, where $z\in\M$ and $\chi\in\U$, instead of using $\chi$ to denote the state-input information as in~\eqref{eq:hybridplant}. Indeed, we observe that such modeling is equivalent to the approach in~\eqref{eq:hybridplant}. 

Let $(\M, \U, \pi)$ denote the relaxed $C^1$-control bundle on which a continuous-time plant evolves. Its state is denoted by $z\in \M$ and the state-dependent input by $\chi \in \U$. Let the continuous-time dynamics be given by the system
\begin{align}
    \label{eq:CT-productDynamics}
    \dot{z} \in \widetilde{F}_P(z, \chi) \qquad (z, \chi) \in \widetilde{C}_P, \quad z = \pi(\chi),
\end{align}
where $\widetilde{C}_P \subset \M \times \U$ denotes the flow set and $\widetilde{F}_P : \M\times \U \setto \T{}{\M}$ denotes the flow map, which satisfies $\widetilde{F}_P(z, \chi) \in \T{z}{\M}$ for each $(z,\chi)\in \M\times \U$. The additional constraint $z = \pi(\chi)$ in~\eqref{eq:CT-productDynamics} ensures that the chosen $\chi$ belongs to the fiber $\pi^{-1}(z)$. Because of this constraint, each $\chi\in\U$ uniquely determines its associated state $z$ and, consequently, the explicit dependence of $\widetilde{F}_P$ on the state $z$ can be equivalently captured using $\chi$ via the relation $z = \pi(\chi)$. Due to this fact, the system~\eqref{eq:CT-productDynamics} can be written equivalently as the constrained continuous-time plant
\begin{align}
    \label{eq:Sigma_C}
    \Sigma_C \;\;:\;\; \dot{z} \in F_P(\chi)\qquad \chi \in C_P,
\end{align}
where $z\in \M$ denotes the state, $\chi \in \pi^{-1}(z) \subset \U$ denotes the state-dependent input, $C_P \subset \U$ denotes the flow set, and $F_P : \U \setto \T{}{\M}$ denotes the flow map that satisfies $F_P(\chi) \subset \T{\pi(\chi)}{\M}$ for each $\chi\in \U$. Note that $F_P$ does not explicitly depend on the state $z$. Therefore, the notation used for $\Sigma_C$ suppresses the redundant state argument in the flow map while describing precisely the same flow dynamics as the explicit $(z,\chi)$ representation in~\eqref{eq:CT-productDynamics}.

\begin{example}[Nontrivial bundle]
    \label{ex:integrator}
    Consider the unconstrained integrator dynamics on a $C^1$-manifold $\M$ as
    \begin{align}
    \label{eq:integrator}
        \dot{z} = u,
    \end{align}
    where ${z}\in \M$ denotes the state and $u \in \T{{z}}{\M}$ denotes the state-dependent control input. In particular, $u \in \T{}{\M}$ and ${z} = \pi_{\M}(u)$, where $\pi_{\M} : \T{}{\M} \to \M$ denotes the canonical tangent bundle projection map. Therefore, with $F_P : \T{}{\M} \to \T{}{\M}$ as the identity map and $C_P \coloneqq \T{}{\M}$,~\eqref{eq:integrator} is modeled in the form $\Sigma_C$ as in~\eqref{eq:Sigma_C} on the relaxed $C^1$-control bundle $(\M, \T{}{\M}, \pi_{\M})$, where $\T{}{\M}$ is a topological space as it is a manifold~\cite[Proposition~3.18]{Lee} and $\pi_{\M}$ is a surjection by definition. Note that, for each ${z}\in \M$, the fiber over ${z}$, namely, $\pi_{\M}^{-1}({z})$, is the tangent space $\T{{z}}{\M}$. 
\end{example}

\begin{example}[Trivial bundle]
\label{ex:attitude-kinematics}
    Consider the set of three-dimensional rotation matrices $\SOthree \coloneqq \{R \in \R{3\times 3} : RR^\top = R^\top R= I, \operatorname{det} R = 1\}$. Then, the attitude kinematics is given by
    \begin{align}
    \label{eq:attitutde-kinematics}
        \dot{R} = R \Omega_{\times}
    \end{align}
    where $R \in \M \coloneqq \SOthree$ denotes the state, $\Omega \in \R{3}$ denotes the control input, and $\times$ is the cross product operator such that, for each $v, w\in \R{3}$, $v_{\times}\in \R{3\times 3}$ such that $v_\times w = v\times w$. Then, with $\V\coloneqq \SOthree\times \R{3}$, ~\eqref{eq:attitutde-kinematics} is modeled on the relaxed $C^1$-control bundle $(\SOthree, \V, \operatorname{pr}_1)$, where $\operatorname{pr}_1 : \V \to \SOthree$ is the projection map defined as $\operatorname{pr}_1(R,\Omega) = R$ for each $(R, \Omega) \in \V$. Note that, as $\V$ is a product space, it is a \emph{trivial bundle}~\cite[p.250]{Lee}. Namely, the fiber $\operatorname{pr}_1^{-1}(R)$ over each point $R \in \M$ remains constant equal to $\R{3}$. As a result, points in $\V$ can be denoted as tuples $(R, \Omega)\in \V$. Such tuple-based representation of state-input pairs mimics the Euclidean modeling of dynamical systems~\cite{Sanfelice_HybridFeedbackControl_2021}, and was unavailable in Example~\ref{ex:integrator} due to lack of a trivial bundle structure.
\end{example}

Similar to the continuous-time plant $\Sigma_C$, a discrete-time plant $\Sigma_D$ is modeled on a relaxed $C^1$-control bundle $(\M, \U, \pi)$ via a set-valued map $G_P : \U \setto {\M}$ and a constraint set $D_P \subset \U$. This system, with state ${z}\in \M$ and input $\chi \in \pi^{-1}({z}) \subset \U$, is denoted by
\begin{align}
    \label{eq:Sigma_D}
    \Sigma_D \;\;:\;\;  {z}^+ \in G_P(\chi) \quad \chi \in D_P.
\end{align}

Examples~\ref{ex:integrator} and~\ref{ex:attitude-kinematics} can be similarly extended to model discrete-time systems with inputs on manifolds.

\end{appendices}

\bibliography{sn-bibliography}
\end{document}